\documentclass{siamltex}
\usepackage{graphicx}
\usepackage{float}
\usepackage{amsmath}
\usepackage{amsfonts}
\usepackage{algorithm}
\usepackage{cite}
\usepackage{todonotes}
\usepackage{siunitx}
\usepackage{amssymb}
\usepackage{tikz}
\usepackage[colorlinks,linkcolor=blue,urlcolor=blue,citecolor=blue]{hyperref}

\newcommand{\R}{{\mathbb R}}
\newcommand{\dealii}{{\sc deal.II}}

\newcommand{\IAC}{I\!AC}

\newtheorem{remark}[theorem]{Remark}

\newcounter{row}
\newcounter{col}

\title{A benchmark for the Bayesian inversion of coefficients in
  partial differential equations}

\author{David Aristoff\thanks{Department of Mathematics, 1874 Campus Delivery, Colorado
    State University; Fort Collins, CO 80523; USA (\texttt{aristoff@math.colostate.edu}).}
        \and
        Wolfgang Bangerth\thanks{Department of Mathematics, Department of Geosciences, 1874 Campus Delivery, Colorado
    State University; Fort Collins, CO 80523; USA (\texttt{bangerth@colostate.edu}).}}

\begin{document}

\maketitle

\begin{abstract}
  Bayesian methods have been widely used in the last two decades to
  infer statistical properties of spatially variable coefficients in partial differential
  equations from measurements of the solutions of these
  equations. Yet, in many cases the number of variables used to
  parameterize these coefficients is large, and obtaining meaningful
  statistics of their values is difficult using simple sampling
  methods such as the basic Metropolis-Hastings (MH) algorithm -- in particular if
  the inverse problem is ill-conditioned or ill-posed. As a
  consequence, many advanced sampling methods have been described in
  the literature that converge faster than MH, for example by
  exploiting hierarchies of statistical models or hierarchies of discretizations of the underlying
  differential equation.

  At the same time, it remains difficult for the reader of the
  literature to quantify the advantages of these algorithms because
  there is no commonly used benchmark. This paper presents a benchmark
  Bayesian inverse problem -- namely, the
  determination of a spatially-variable coefficient, discretized by 64 values, in a Poisson equation, based on point
  measurements of the solution -- that fills the gap
  between widely used simple test cases (such as superpositions of
  Gaussians) and real applications that are difficult to replicate
  for developers of sampling algorithms. We provide a complete
  description of the test case, and provide an open source
  implementation that can serve as the basis for further
  experiments. We have also computed $2\times 10^{11}$ samples, at a
  cost of some 30 CPU years, of the posterior probability distribution
  from which we have generated detailed and accurate statistics against which other
  sampling algorithms can be tested.
\end{abstract}

\begin{AMS}
  65N21, 35R30, 74G75
\end{AMS}

\section{Introduction}

Inverse problems are parameter estimation problems in which one wants
to determine unknown, spatially variable material parameters in a partial differential equation
(PDE) based on measurements of the solution. In the deterministic approach,
one in essence seeks that set of parameters for which the solution of the
PDE would best match the measured values; this approach is widely used
in many applications. On the other hand, the Bayesian approach to inverse
problems recognizes that all measurements are subject to measurement
errors and that models are also inexact; as a consequence, we ought
to pose the inverse problem as one that seeks a probability distribution
describing how likely it is that parameter values lie in a given interval
or set. This generalization of the perspective on inverse problems has
long roots, but first came to the attention of the wider scientific
community through a 1987 book by Tarantola~\cite{Tarantola87}. It was later
followed by a significantly revised and more accessible version by the same
author \cite{Tarantola2005} as well as numerous other books on the subject;
we mention \cite{Kaipio2005} as one example, along with \cite{Allmaras2013,Aguilar2015,dashti2015bayesian} for tutorial-style introductions to the topic. The Bayesian approach to inverse problems has been used in a wide variety of inverse applications: Too many to mention in detail, but including 
acoustics~\cite{BuiThanh2012}, flow in the Earth's mantle~\cite{Worthen2014}, laminar and turbulent flow~\cite{chen2019}, ice sheet modeling~\cite{Petra2014}, astronomy~\cite{astro}, chemistry~\cite{friedman2000using,neal1993bayesian}, and groundwater modeling~\cite{Jiang2006}.

From a computational perspective,
the primary challenge in Bayesian inverse problems is that
after discretizing the spatially variable parameters one seeks to
infer, one generally ends up with trying to characterize a finite- but
high-dimensional probability distribution $\pi(\theta)$ that describes
the relative likelihood of parameters $\theta$. In particular, we are
typically interested in computing the mean and standard deviation of this
probability distribution (i.e., which set of parameters $\left<\theta\right>$
\textit{on average} fits the measured data best, and what we know about
its variability given the uncertainty in the measured data). Computing
these integral quantities in high dimensional spaces can only be done
through sampling methods such as Monte Carlo Markov Chain (MCMC) algorithms.
On the other hand, sampling in high-dimensional spaces often suffers from a number of problems: (i) Long burn-in times until a chain finally finds the
region where the probability distribution $\pi(\theta)$ has values
substantially different from zero; (ii) long autocorrelation length scales
if $\pi(\theta)$ represents elongated, curved, or ``ridged'' distributions; (iii) for some inverse problems, multi-modality of $\pi(\theta)$ is also a problem that complicates the interpretation of the posterior probability distribution; (iv) in many \textit{ill-posed} inverse problems, parameters  have large variances that result in rather slow convergence to reliable and accurate means.

In most high-dimensional applications, the result of these issues is then
that one needs very large numbers of samples to accurately characterize
the desired probability distribution. In the context of inverse problems,
this problem is compounded by the fact that the generation of every sample
requires the solution of the forward problem -- i.e., generally, the
expensive numerical solution of a PDE. As a consequence, the solution
of Bayesian inverse problems is computationally exceptionally expensive.

The community has stepped up to this challenge over the past two decades:
Numerous algorithms have been developed to make the sampling process
more efficient. Starting from the simplest sampler, the
Metropolis-Hastings algorithms with a symmetric proposal distribution
\cite{hastings1970monte}, ideas to alleviate some of the problems include nonsymmetric
proposal distributions~\cite{roberts1996exponential}, delayed rejection~\cite{tierney1999some}, non-reversible samplers~\cite{duncan2016variance}, piecewise deterministic Markov processes~\cite{bierkens2019zig,vanetti2017piecewise} including Hamiltonian Monte Carlo~\cite{neal2011mcmc}, adaptive methods~\cite{Haario_adaptive,atchade2005adaptive,roberts2009examples}, 
randomize-then-optimize methods \cite{Bardsley2014,Bardsley2015,Bardlsey2018}, 
affine invariant samplers~\cite{goodman2010ensemble,foreman2013emcee},
and combinations thereof~\cite{darve2008adaptive,haario2006dram}.

Other approaches introduce parallelism (e.g., the differential evolution
method and variations \cite{Braak2006,terBraak2008,Vrugt2009}), or hierarchies of models (see,
for example, the survey by Peherstorfer et al. \cite{Peherstorfer2018} and references therein, as well as \cite{Seo2020,Fleeter2020,Christen2005,Watzenig_2009,Fox1997}). 
Yet other methods exploit the fact
that discretizing the underlying PDE gives rise to a natural
multilevel hierarchy (see \cite{Dodwell2015} among many others) or that one can use the structure of the discretized PDE for efficient sampling algorithms \cite{hippylib,Martin2012}.

Many of these methods are likely vastly faster than the simplest sampling
methods that are often used. Yet, the availability of a whole zoo of possible
methods and their various possible combinations has also made it difficult to
assess which method really should be used if one wanted to solve a
particular inverse problem, and there is no consensus in
the community on this topic. Underlying this lack of consensus is that
\textit{there is no widely used benchmark} for Bayesian inverse problems:
Most of the papers above demonstrate the qualities of their particular
innovation using some small but artificial test cases such as a
superposition of Gaussians, and often a more elaborate application that
is insufficiently well-described and often too complex for others to
reproduce. As a consequence, the literature contains few examples of
comparisons of algorithms \textit{using test cases that reflect the
properties of actual inverse problems}.

Our contribution here seeks to address this lack of widely used benchmarks.
In particular:
\begin{itemize}
    \item We provide a complete description of a benchmark that involves
      characterizing a posterior probability distribution $\pi(\theta)$
      on a 64-dimensional parameter space that results from inverting
      data for a discretized coefficient in a Poisson equation.
    \item We explain in detail why this benchmark is at once simple enough
      to make reproduction by others possible, yet difficult enough
      to reflect the real challenges one faces when solving Bayesian
      inverse problems.
    \item We provide highly accurate statistics for $\pi(\theta)$
      that allow others to assess the correctness of their own algorithms
      and implementations. We also provide a performance profile
      for a simple Metropolis-Hastings sampler as a baseline against
      which other methods can be compared.
\end{itemize}
To make adoption of this benchmark simpler, we also provide an open-source
implementation of the benchmark that can be adapted to experimentation
on other sampling methods with relative ease.

The remainder of this paper is structured as follows: In Section~\ref{sec:benchmark-description}, we provide a complete description of the benchmark. In Section~\ref{sec:assessment}, we then evaluate
highly accurate statistics of the probability distribution that solves
the benchmark, based on \num{2e11} samples we have computed.
Section~\ref{sec:hope} provides a short discussion of what we hope
the benchmark will achieve, along with our conclusions.
An appendix presents details of our implementation of the benchmark (Appendix~\ref{sec:software}), discusses a simple 1d benchmark for which one can find solutions in a much cheaper way (Appendix~\ref{sec:1d}), and provides some statistical background relevant to Section~\ref{sec:assessment} (in Appendix~\ref{sec:background}).

\section{The benchmark for sampling algorithms for inverse problems}
\label{sec:benchmark-description}

\subsection{Design criteria}
In the design of the benchmark described in this paper, we were guided
by the following principle:
\begin{quote}
  \it
  A good benchmark is neither too simple nor too complicated. It also
  needs to reflect properties of real-world applications.
\end{quote}
Specifically, we design a benchmark for inferring posterior
probability distributions using sampling algorithms that correspond to
the Bayesian inversion of coefficients in partial differential
equations. In other words, cases where the relative posterior
likelihood is computed by comparing (functionals of) the forward
solution of partial differential equations with (simulations of)
measured data.

The literature has many examples of papers that consider sampling
algorithms for such problems (see the references in the introduction).
However, they are typically tested only
on cases that fall into the following two categories:
\begin{itemize}
  \item Simple posterior probability density functions (PDFs) that are
    given by explicitly known expressions such as Gaussians or
    superpositions of Gaussians. There are advantages to such test
    cases: (i) the probability distributions are cheap to evaluate,
    and it is consequently possible to create essentially unlimited
    numbers of samples; (ii) because the PDF is exactly known, exact
    values for statistics such as the mean, covariances, or maximum
    likelihood (MAP) points are often computable exactly, facilitating the
    quantitative assessment of convergence of sampling schemes. On the
    other hand, these test cases are often so simple that \textit{any}
    reasonable sampling algorithm converges relatively quickly, making
    true comparisons between different algorithms difficult. More
    importantly, however, such simple test cases \textit{do not
      reflect real-world properties} of inverse problems: Most inverse
    problems are ill-posed, nonlinear, and high-dimensional. They are
    often unimodal, but with PDFs that are typically quite insensitive
    along certain directions in parameter space, reflecting the
    ill-posedness of the underlying problem. Because real-world
    problems are so different from simple artificial test cases, it is
    difficult to draw conclusions from the performance of a new
    sampling algorithm when applied to a simple test case.

  \item Complex applications, such as the determination of the
    spatially variable oil reservoir permeability from the production
    history of an oil field, or the determination of seismic wave
    speeds from travel times of earthquake waves from their source to
    receivers. Such applications are of course the target for applying
    advanced sampling methods, but they make for poor benchmarks
    because they are very difficult to replicate by other authors. As
    a consequence, they are almost exclusively used in only the single
    paper in which a new sampling algorithm is first described, and
    it is difficult for others to compare this new sampling
    algorithm against previous ones, since they have not
    been tested against the same, realistic benchmark.
\end{itemize}

We position the benchmark in this paper between these
extremes. Specifically, we have set out to achieve the following goals:
\begin{itemize}
\item \textit{Reflect real properties of inverse problems:} Our
  benchmark should reflect properties one would expect from real
  applications such as the permeability or seismic wave speed
  determinations mentioned above. We do not really know what these
  properties are, but intuition and knowledge of the literature
  suggest that they include very elongated and nonlinear probability
  distributions, quite unlike Gaussians or their superpositions. 
  In order for our benchmark to
  reflect these properties, we base it on a
  partial differential equation.
\item \textit{High dimensional:} Inverse problems are originally
  infinite dimensional, i.e., we seek parameters that are functions of
  space and/or time. In practice, these need to be discretized,
  leading to finite- but often high-dimensional
  problems. It is well understood that the resulting curse of
  dimensionality leads to practical problems that often make the Bayesian
  inverse problem extremely expensive to solve. At the same time, we
  want to reflect these difficulties in our benchmark.
\item \textit{Computable with acceptable effort:} A benchmark needs to
  have a solution that is known to an accuracy sufficiently good to
  compare against. This implies that it can't be so expensive that we
  can only compute a few thousand or tens of thousands of samples of
  the posterior probability distribution. This rules out most real
  applications for which each forward solution, even on parallel
  computers may take minutes, hours, or even longer. Rather, we need a
  problem that can be solved in at most a second on a single processor
  to allow the generation of a substantial number of samples.
\item \textit{Reproducible:} To be usable by anyone, a benchmark needs
  to be completely specified in all of its details. It also needs to
  be simple enough so that others can implement it with reasonable effort.
\item \textit{Available:} An important component of this paper is that
  we make the software that implements the benchmark available as open
  source, see Appendix~\ref{sec:software}. In particular, the code
  is written in a modular way that allows evaluating the posterior
  probability density for a given set of parameter values -- i.e., the key
  operation of all sampling methods. The code is also written in such
  a way that it is easy to use in a multilevel sampling scheme where
  the forward problem is solved with a hierarchy of successively more
  accurate approximations.
\end{itemize}

\subsection{Description of the benchmark}

Given the design criteria discussed in the previous subsection, let us
now present the details of the benchmark. Specifically, we seek
(statistics of) a non-normalized
posterior probability distribution $\pi(\theta|\hat z)$ on a parameter space
$\theta\in\Theta=\R^{64}$ of modestly high dimension 64 -- large
enough to be interesting, while small enough to remain feasible for a benchmark. Here, we
think of $\hat z$ as a set of measurements made on a physical system that
is used to infer information about the internal parameters $\theta$ of
the system. As is common
in Bayesian inverse problems, $\pi(\theta|\hat z)$ is defined as the product
of a likelihood times a prior probability:
\begin{align}
  \label{eq:posterior}
  \pi(\theta|\hat z) \propto L(\hat z|\theta) \, \pi_\text{pr}(\theta).
\end{align}
Here, $L(z|\theta)$ describes how likely it would be to measure values
$z$ if $\theta$ were the ``true'' values of the internal
parameters. $\pi_\text{pr}$ is a (not necessarily normalized) probability
distribution encoding our prior beliefs about the parameters.
A complete description of the benchmark then requires us to describe
the values of $z$ and ways to evaluate the functions $L$ and
$\pi_\text{pr}$. We will split the definition of $L$ into a discussion
of the forward model and a statistical model of measurements in the
following.

\subsubsection{The forward model}
\label{sec:forward-model}
\begin{figure}
  \centering
  
  \phantom{.}
  \hfill
  \begin{tikzpicture}[scale=0.65]
    \begin{scope}
      \draw (0, 0) grid (8, 8);

      \tiny

      \begin{scope}[gray]
        \foreach \i in {1, 2, 3, 4, 5, 6, 7, 8} {
          \foreach \j in {1, 2, 3, 4, 5, 6, 7, 8} {
            \edef\x{\i - 0.5}
            \edef\y{\j - 0.5}
            \pgfmathtruncatemacro\n{(\j-1) * 8 + (\i-1)}
            \node[anchor=center] at (\x, \y) {\n};
          }
        }
      \end{scope}
    \end{scope}
  \end{tikzpicture}
  \hfill
  \begin{tikzpicture}[scale=0.65]
    \begin{scope}
      \draw (0, 0) grid (8, 8);

      \tiny

      \begin{scope}[blue]
        \foreach \i in {1, 2, 3, 4, 5, 6, 7, 8, 9, 10, 11, 12, 13} {
          \foreach \j in {1, 2, 3, 4, 5, 6, 7, 8, 9, 10, 11, 12, 13} {
            \edef\x{\i * 8.0 / 14.0}
            \edef\y{\j * 8.0 / 14.0}
            \edef\n{\i}
            \node[anchor=center] at (\x, \y) {$\bullet$};
          }
        }
      \end{scope}

      \begin{scope}[gray]
        \foreach \i in {1, 2, 3, 4, 5, 6, 7, 8, 9, 10, 11, 12, 13} {
          \foreach \j in {1, 2, 3, 4, 5, 6, 7, 8, 9, 10, 11, 12, 13} {
            \edef\x{\i * 8.0 / 14.0}
            \edef\y{\j * 8.0 / 14.0}
            \pgfmathtruncatemacro\n{(\j-1) * 13 + (\i-1)}
            \node[anchor=west] at (\x, \y) {\scalebox{0.75}{$\!\small\n$}};
          }
        }
      \end{scope}
      
    \end{scope}
  \end{tikzpicture}
  \hfill
  \phantom{.}

  \caption{\it Left: Numbering of the 64 cells on which the
    parameters are defined. Right: Numbering and locations
    of the $13^2=169$
    evaluation points at which the solution is evaluated.}
  \label{fig:forward}
\end{figure}
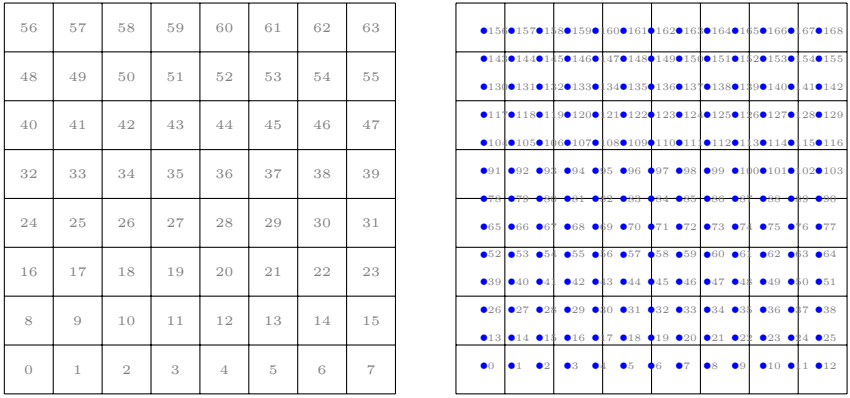

The setting we want to pursue is as follows: Let us imagine a membrane
stretched over a frame that bounds a domain $\Omega$ which, for
simplicity we assume to be the unit square $\Omega=(0,1)^2$. The
membrane is subject to an external vertical force $f(\mathbf x)$ which
for the purpose of this benchmark we choose constant as $f(\mathbf
x)=10$.  Furthermore, the membrane has a spatially variable resistance
$a(\mathbf x)$ to deflection (for example, it may have a variable
thickness or may be made from different materials). In this benchmark,
we assume that $a(\mathbf x)$ is piecewise constant on a uniform $8\times 8$
grid as shown in Fig.~\ref{fig:forward}, with the 64 values that parameterize
$a(\mathbf x)$ given by the elements of the vector
$\theta_0,\ldots,\theta_{63}$ as also indicated in the figure. In
other words, there is a 1:1 relationship between the vector $\theta$
and the piecewise constant coefficient function $a(\mathbf x)=a^\theta(\mathbf x)$.

Then, an
appropriate model to describe the vertical deflection $u(\mathbf x)$
of the membrane would express $u$ as the solution of the following partial
differential equation that generalizes the Poisson equation:
\begin{align}
  \label{eq:laplace-1}
  -\nabla \cdot \left[ a(\mathbf x) \nabla u(\mathbf x) \right]
  &= f(\mathbf x)
  &\qquad\qquad &\text{in $\Omega$},
  \\
  \label{eq:laplace-2}
  u(\mathbf x)
  &= 0
  &\qquad\qquad &\text{on $\partial\Omega$}.
\end{align}
This model is of course not exactly solvable. But its solution can be
approximated using discretization. The way we define the likelihood
$L$ then requires us to specify exactly how we discretize this
model. Concretely, we define $u_h(\mathbf x)$ as the solution of a
finite element discretization of
\eqref{eq:laplace-1}--\eqref{eq:laplace-2} using a uniform $32\times
32$ mesh and a $Q_1$ (bilinear) element. Because $f(\mathbf x)$ is
given, and because there is a 1:1 relationship between $\theta$ and
$a(\mathbf x)$, this discretized model then implies that for each
$\theta$ we can find a $u_h(\mathbf x)=u_h^\theta(\mathbf x)$ that can
be thought of as being parameterized using the \num{1089} degrees of
freedom of the $Q_1$ discretization on the $32\times 32$ mesh. (However, of
these \num{1089} degrees of freedom, 128 are on the boundary and are
constrained to zero.) In other words, using the $Q_1$ shape functions
$\varphi_k(\mathbf x)$,
we can express $u_h(\mathbf x)=\sum_{k=0}^{1088} U_k \varphi_k(\mathbf
x)$. It is important to stress that the mapping
$\theta \mapsto u_h^\theta(\mathbf x)$ (or equivalently,
$\theta\mapsto U^\theta$) is \textit{nonlinear}. 

The function $u_h^\theta(\mathbf x)$ can be thought of as the predicted
displacement at every point $\mathbf x\in\Omega$ if $\theta$
represented the spatially variable stiffness coefficient of the
membrane. In practice, however, we can only measure finitely many
things, and consequently define a measurement operator ${\cal M}: u_h \mapsto
z\in \R^{169}$ that evaluates $u_h$ on a uniform $13\times 13$ grid of
points $\mathbf x_k \in \Omega$ so that $\mathbf
x_k=\left(\frac{i}{13+1},\frac{j}{13+1}\right), 1\le i,j\le 13$
with $k=13(i-1)+(j-1)$.
The
locations of these points are also indicated in Fig.~\ref{fig:forward}. This
last step then defines a linear mapping. Because of the
equivalence between the function $u_h$ and its nodal vector $U$, the
linearity of the measurement operator implies that we can write $z=MU$
with a matrix $M\in \R^{169\times 1089}$ that is given by
$M_{kl}=\varphi_l(\mathbf x_k)$.

In summary, a parameter vector $\theta\in \R^{64}$ then predicts measurements
$z^\theta \in \R^{169}$ using the following chain of maps:
\begin{align}
  \label{eq:forward-model}
  \theta \mapsto a^\theta(\mathbf x) \mapsto U^\theta
  \mapsto z^\theta.
\end{align}
The mapping $\theta \mapsto z^\theta$ is commonly called the ``forward
model'' as it \textit{predicts} measurements $z^\theta$ if we knew the
parameter values $\theta$. The ``inverse problem'' is then of course
the inverse operation: to \textit{infer} the parameters $\theta$ that
describe a system based on measurements $z$ of its state $u$.

All of the steps of the forward model have been precisely defined
above and are easily computable with some basic knowledge of finite
element methods (or using the code discussed in
Appendix~\ref{sec:software}). The expensive step is to solve for
the nodal vector $U^\theta$, as this requires the assembly and
solution of a linear system of size \num{1089}.

\begin{remark}
  The $32\times 32$ mesh to define the forward model is chosen
  sufficiently fine to resolve the exact solution $u$ reasonably
  well. At the same time, it is coarse enough to allow for the rapid
  evaluation of the solution -- even a rather simple implementation
  should yield a solution in less than a second, and a highly
  optimized implementation such as the one discussed in
  Appendix~\ref{sec:software-forward} will be able to do so in less than 5 milliseconds
  on modern hardware. As a consequence, this choice of mesh allows for
  computing a large number of samples, and consequently accurate
  quantitative comparisons of sampling algorithms.

  We also mention that the $32\times 32$ mesh for $u_h(\mathbf x)$ is
  twice more globally refined than the $8\times 8$ mesh used to define
  $a(\mathbf x)$ in terms of $\theta$. It is clear to practitioners of
  finite element discretizations of partial differential equations
  that the mesh for $u_h$ must be at the very least as fine as the one
  for the coefficient $a^\theta$ to obtain any kind of accuracy. On the other hand, these
  choices then leave room for a hierarchy of models in which the
  forward model uses $8\times 8$, $16\times 16$, and $32\times 32$
  meshes; we expect that multilevel sampling methods will use this hierarchy to good effect.
\end{remark}

\begin{remark}
  \label{remark:size}
  In our experiments, we will choose the values of $\theta$ (and
  consequently of $a(\mathbf x)$) clustered around one. With the
  choice $f=10$ mentioned above, this leads to a solution $u(\mathbf
  x)$ with values in the range $0\ldots 0.95$. This then also implies
  that we should think of the numerical magnitude of our measurements
  $z^\theta_k$ as ${\cal O}(1)$.
\end{remark}

\subsubsection{The likelihood $L(z|\theta)$}
Given the predicted measurements $z^\theta$ that correspond to a
given set of parameters $\theta$, the likelihood $L(z|\theta)$ can
be thought of as expressing the (non-normalized) probability of
actually obtaining $z$ in a measurement if $\theta$ were the ``correct'' set of parameters. This is a statement that
encodes the measurement error of our measurement device.

For the purposes of this benchmark, we assume that these measurement
errors are identical and independently distributed for all 169 measurement
points. More specifically, we define the likelihood as the
following (non-normalized) probability function:
\begin{align}
  \label{eq:likelihood}
  L(z|\theta)
  =
  \exp\left(-\frac{\|z-z^\theta\|^2}{2\sigma^2}\right)
  =
  \prod_{k=0}^{168} \exp\left(-\frac{(z_k-z^\theta_k)^2}{2\sigma^2}\right),
\end{align}
where we set $\sigma=0.05$ and where $z^\theta$ is related to $\theta$
using the chain \eqref{eq:forward-model}.

\begin{remark}
  \label{remark:sigma}
  We can think of \eqref{eq:likelihood} as encoding our belief that
  our measurement system produces a Gaussian-distributed measurement $z_k \sim
  N(z^\theta_k, \sigma)$. Given that 
  $z^\theta_k={\cal O}(1)$, $\sigma=0.05$ implies a
  measurement error of 5\%. This is clearly much larger than the
  accuracy with which one would be able to determine the deflection of
  a membrane in practice. On the other hand, we have chosen $\sigma$
  this large to ensure that the Bayesian inverse problem does not lead
  to a probability distribution $\pi(\theta|\hat z)$ that is so
  narrowly centered around a value $\bar\theta$ that the mapping
  $\theta\mapsto z^\theta$ can be linearized around $\bar\theta$ -- in
  which case the likelihood $L(\hat z|\theta)$ would
  become Gaussian, as also discussed in Appendix~\ref{sec:1d}. We will demonstrate in Section~\ref{sec:nongaussianty} that
  indeed $\pi(\theta|\hat z)$ is not Gaussian and, moreover, is large
  along a curved ridge that can not easily be approximated by a
  Gaussian either.
\end{remark}

\subsubsection{The prior probability $\pi_\text{pr}(\theta)$}
\label{sec:prior}
Our next task is to describe our prior beliefs for the values of the
parameters. Given that the 64 values of $\theta$ describe the
stiffness coefficient of a membrane, it is clear that they must be
positive. Furthermore, as with many mechanical properties that can have values over vast ranges,%
\footnote{For example, the Young's modulus that is related to the stiffness of a membrane, can range from 0.01 GPa for rubber to 200 GPa for typical steels. Similarly, the permeability of typical
oil reservoir rocks can range from 1 to 1000 millidarcies.}
reasonable priors are typically posed on the ``order of magnitude'' (that is, the \textit{logarithm}), not the size of the coefficient itself.
We express this through the
following (non-normalized) probability distribution:
\begin{align}
  \label{eq:prior}
  \pi_\text{pr}(\theta) = \prod_{i=0}^{63} \exp\left(-\frac{(\ln(\theta_i)-\ln(1))^2}{2\sigma_\text{pr}^2}\right),
\end{align}
where we choose $\sigma_\text{pr}=2$. We
recognize the prior density of $\ln(\theta_k)$ as a Gaussian with mean $\sigma_\text{pr}^2$ and 
standard deviation $\sigma_\text{pr}$. 

Because this prior distribution is posed on the logarithm of the parameters, the prior on the parameters themselves is very heavy-tailed, with mean values $\left<\theta_k\right>_{\pi_\text{pr}}$ for each component much larger than the value at which $\pi_\text{pr}$ takes on its maximum (which is at $\theta_k=1$): Indeed, the
mean of each $\theta_k$ with respect to $\pi_{\textup{pr}}$ is about $403.43$.

We note that this prior probability is quite weak and, in particular, does not assume any (spatial) correlation between parameters as is often done in inverse problems \cite{hippylib,Stuart2010,Kaipio2005}. The correlations we will observe in our posterior probability (see Section~\ref{sec:covariance}) are therefore a consequence of the likelihood function only.

\subsubsection{The ``true'' measurements $\hat z$}
\label{sec:true-measurements}
The last piece necessary to describe the complete benchmark is the
choice of the ``true'' measurements $\hat z$ that we want to use to
infer the statistical properties of the parameters $\theta$.
For the purposes of this benchmark, we will use the 169 values for $\hat
z$ given in Table~\ref{table:hat-z}.

\begin{table}
  \caption{\it The ``true'' measurement values $\hat z_k, k=0,\ldots,
    168$ used in the benchmark. The values are also available in the
    electronic supplemental material and are shown in full double
    precision accuracy to allow for exact reproduction of the 
    benchmark.}
  \tiny
  \centering
  \begin{tabular}{rl|rl|rl}
    $\hat z_{0}$ & \num{0.06076511762259369} &     $\hat z_{60}$ & \num{0.6235300574156917} &   $\hat z_{120}$ & \num{0.5140091754526943} \\
    & \num{0.09601910120848481} &     & \num{0.5559332704045935} &   & \num{0.5559332704045969} \\
    & \num{0.1238852517838584} &     & \num{0.4670304994474178} &     & \num{0.5677443693743304} \\
    & \num{0.1495184117375201} &     & \num{0.3499809143811} &     & \num{0.5478251665295453} \\
    & \num{0.1841596127549784} &     & \num{0.19688263746294} &     & \num{0.4895759664908982} \\
    & \num{0.2174525028261122} &     & \num{0.2174525028261253} &     & \num{0.4109641743019171} \\
    & \num{0.2250996160898698} &     & \num{0.4122329537843404} &     & \num{0.395727260284338} \\
    & \num{0.2197954769002993} &     & \num{0.5779452419831566} &     & \num{0.3778949322004734} \\
    & \num{0.2074695698370926} &     & \num{0.6859683749254372} &     & \num{0.3596268271857124} \\
    & \num{0.1889996477663016} &     & \num{0.7373108331396063} &     & \num{0.2191250268948948} \\
    $\hat z_{10}$ & \num{0.1632722532153726} &     $\hat z_{70}$ & \num{0.7458811983178246} &     $\hat z_{130}$ & \num{0.1632722532153683} \\
    & \num{0.1276782480038186} &     & \num{0.7278968022406559} &     & \num{0.2850397806663325} \\
    & \num{0.07711845915789312} &     & \num{0.6904793535357751} &     & \num{0.373006008206081} \\
    & \num{0.09601910120848552} &     & \num{0.6369176452710288} &     & \num{0.4325325506354207} \\
    & \num{0.2000589533367983} &     & \num{0.5677443693743215} &     & \num{0.4670304994474315} \\
    & \num{0.3385592591951766} &     & \num{0.4784738764865867} &     & \num{0.4784738764866023} \\
    & \num{0.3934300024647806} &     & \num{0.3602190632823262} &     & \num{0.4677122687599041} \\
    & \num{0.4040223892461541} &     & \num{0.2031792054737325} &     & \num{0.4341716881061055} \\
    & \num{0.4122329537843092} &     & \num{0.2250996160898818} &     & \num{0.388186479011099} \\
    & \num{0.4100480091545554} &     & \num{0.4100480091545787} &     & \num{0.3778949322004602} \\
    $\hat z_{20}$ & \num{0.3949151637189968} &     $\hat z_{80}$ & \num{0.5555615956137137} &     $\hat z_{140}$ & \num{0.3633362567187364} \\
    & \num{0.3697873264791232} &     & \num{0.6561235366960938} &     & \num{0.3464457261905399} \\
    & \num{0.33401826235924} &     & \num{0.7116558878070715} &     & \num{0.2096362321365655} \\
    & \num{0.2850397806663382} &     & \num{0.727896802240657} &     & \num{0.1276782480038148} \\
    & \num{0.2184260032478671} &     & \num{0.7121928678670187} &     & \num{0.2184260032478634} \\
    & \num{0.1271121156350957} &     & \num{0.6712187391428729} &     & \num{0.2821694983395252} \\
    & \num{0.1238852517838611} &     & \num{0.6139157775591492} &     & \num{0.3248315148915535} \\
    & \num{0.3385592591951819} &     & \num{0.5478251665295381} &     & \num{0.3499809143811097} \\
    & \num{0.7119285162766475} &     & \num{0.4677122687599031} &     & \num{0.3602190632823333} \\
    & \num{0.8175712861756428} &     & \num{0.3587654911000848} &     & \num{0.3587654911000799} \\
    $\hat z_{30}$ & \num{0.6836254116578105} &     $\hat z_{90}$ & \num{0.2050734291675918} &     $\hat z_{150}$ & \num{0.3534389974779268} \\
    & \num{0.5779452419831157} &     & \num{0.2197954769003094} &     & \num{0.3642640090182283} \\
    & \num{0.5555615956136897} &     & \num{0.3949151637190157} &     & \num{0.35962682718569} \\
    & \num{0.5285181561736719} &     & \num{0.5285181561736911} &     & \num{0.3464457261905295} \\
    & \num{0.491439702849224} &     & \num{0.6213197201867471} &     & \num{0.3260728953424643} \\
    & \num{0.4409367494853282} &     & \num{0.6745179049094407} &     & \num{0.180670595355394} \\
    & \num{0.3730060082060772} &     & \num{0.690479353535786} &     & \num{0.07711845915789244} \\
    & \num{0.2821694983395214} &     & \num{0.6712187391428787} &     & \num{0.1271121156350963} \\
    & \num{0.1610176733857739} &     & \num{0.6178408289359514} &     & \num{0.1610176733857757} \\
    & \num{0.1495184117375257} &     & \num{0.5453605027237883} &     & \num{0.1834600412730144} \\
    $\hat z_{40}$ & \num{0.3934300024647929} &     $\hat z_{100}$ & \num{0.489575966490909} &     $\hat z_{160}$ & \num{0.1968826374629443} \\
    & \num{0.8175712861756562} &     & \num{0.4341716881061278} &     & \num{0.2031792054737354} \\
    & \num{0.9439154625527653} &     & \num{0.3534389974779456} &     & \num{0.2050734291675885} \\
    & \num{0.8015904115095128} &     & \num{0.2083227496961347} &     & \num{0.2083227496961245} \\
    & \num{0.6859683749254024} &     & \num{0.207469569837099} &     & \num{0.2179599909279998} \\
    & \num{0.6561235366960599} &     & \num{0.3697873264791366} &     & \num{0.2191250268948822} \\
    & \num{0.6213197201867315} &     & \num{0.4914397028492412} &     & \num{0.2096362321365551} \\
    & \num{0.5753611315000049} &     & \num{0.5753611315000203} &     & \num{0.1806705953553887} \\
    & \num{0.5140091754526823} &     & \num{0.6235300574157017} &     $\hat z_{168}$ & \num{0.1067965550010013}\\
    & \num{0.4325325506354165} &     & \num{0.6369176452710497} && \\
    $\hat z_{50}$ & \num{0.3248315148915482} &     $\hat z_{110}$ & \num{0.6139157775591579} && \\
    & \num{0.1834600412730086} &     & \num{0.5453605027237935} && \\
    & \num{0.1841596127549917} &     & \num{0.4336604929612851} && \\
    & \num{0.4040223892461832} &     & \num{0.4109641743019312} && \\
    & \num{0.6836254116578439} &     & \num{0.3881864790111245} && \\
    & \num{0.8015904115095396} &     & \num{0.3642640090182592} && \\
    & \num{0.7870119561144977} &     & \num{0.2179599909280145} && \\
    & \num{0.7373108331395808} &     & \num{0.1889996477663011} && \\
    & \num{0.7116558878070463} &     & \num{0.3340182623592461} && \\
    & \num{0.6745179049094283} &     & \num{0.4409367494853381} &&
  \end{tabular}  \label{table:hat-z}
\end{table}

In some sense, it does not matter where these values come from -- we
could have measured them in an actual experiment, and used these
values to infer the coefficients of the system we measured on. On the
other hand, for the purposes of a benchmark, it might be interesting
to know whether these ``true measurements'' $\hat z$ correspond to a
``true set of parameters'' $\hat\theta$ against which we can compare
statistics such as the mean $\left<\theta\right>$ of the posterior
probability $\pi(\theta|\hat z)$.

\begin{figure}
  \centering
  \begin{tabular}{lp{0.1\textwidth}l}
    \raisebox{0.1\textwidth}{
  \includegraphics[width=0.35\textwidth]{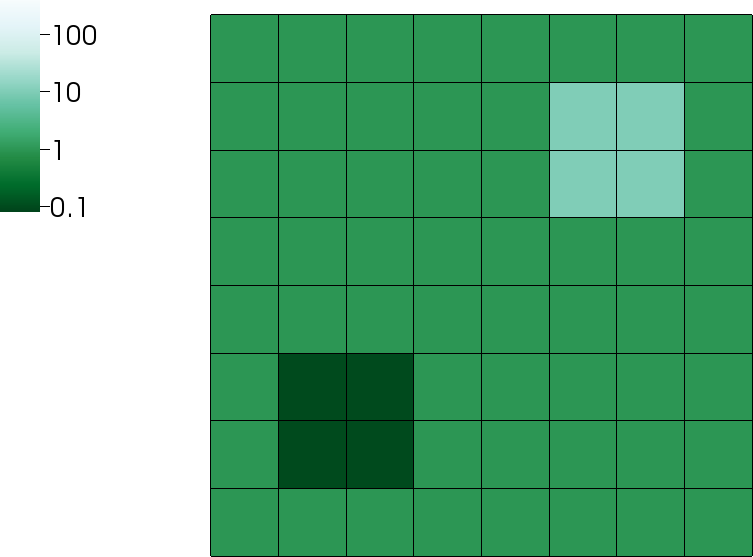}}
  &&
  \includegraphics[width=0.4\textwidth]{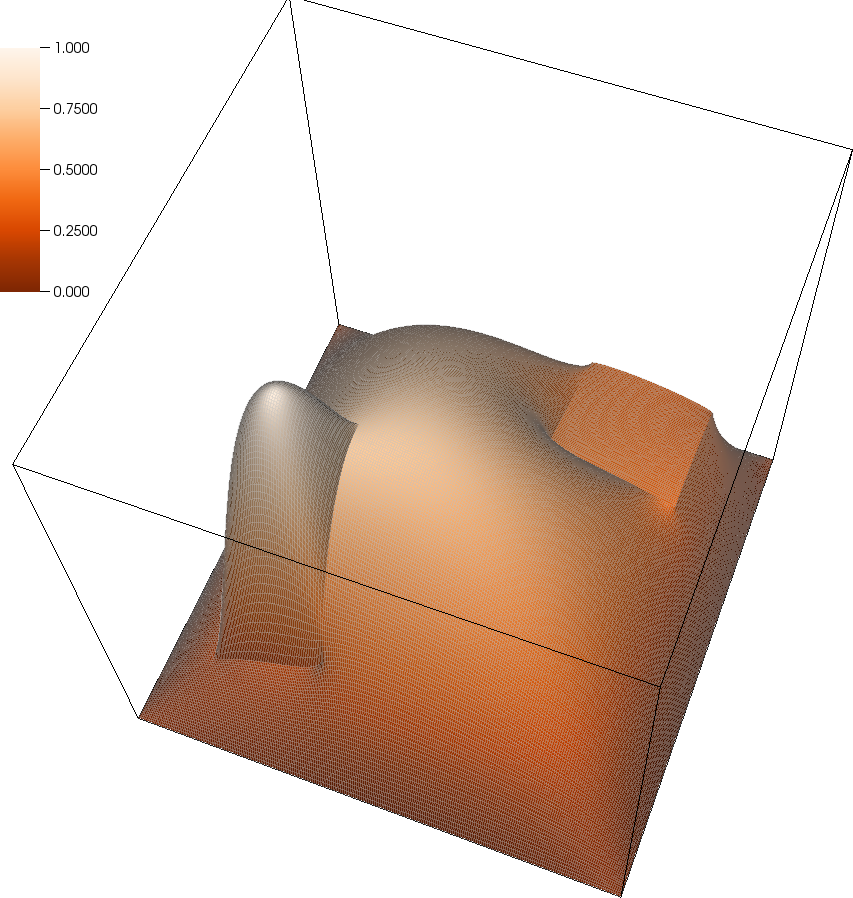}
  \\
  \includegraphics[width=0.4\textwidth]{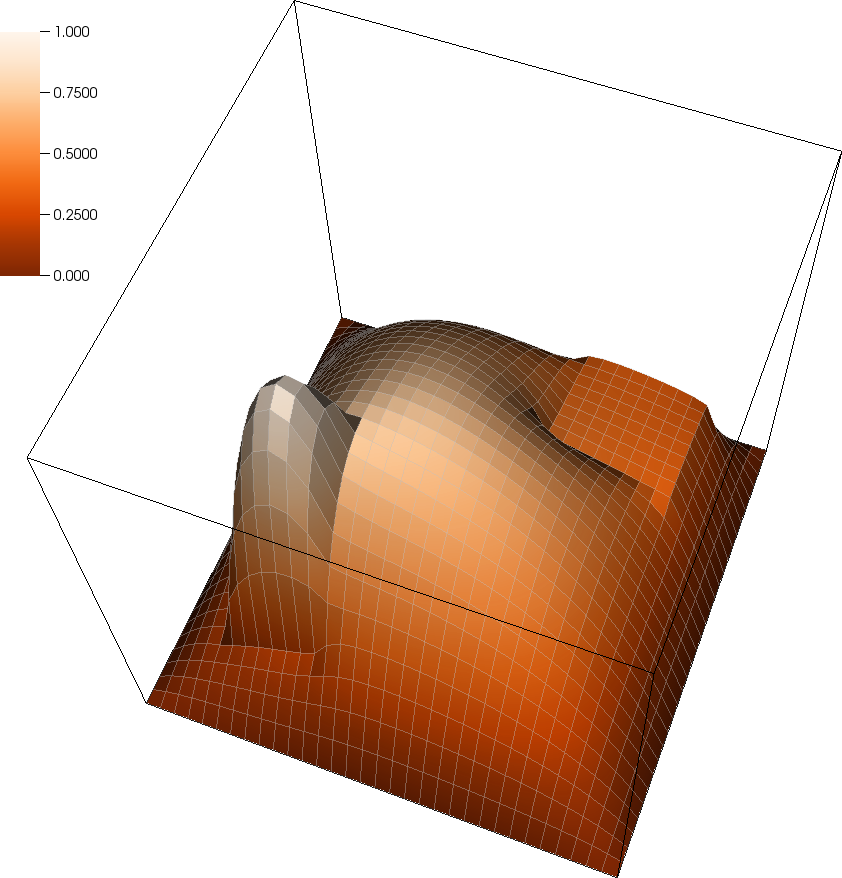}
  &&
  \includegraphics[width=0.4\textwidth]{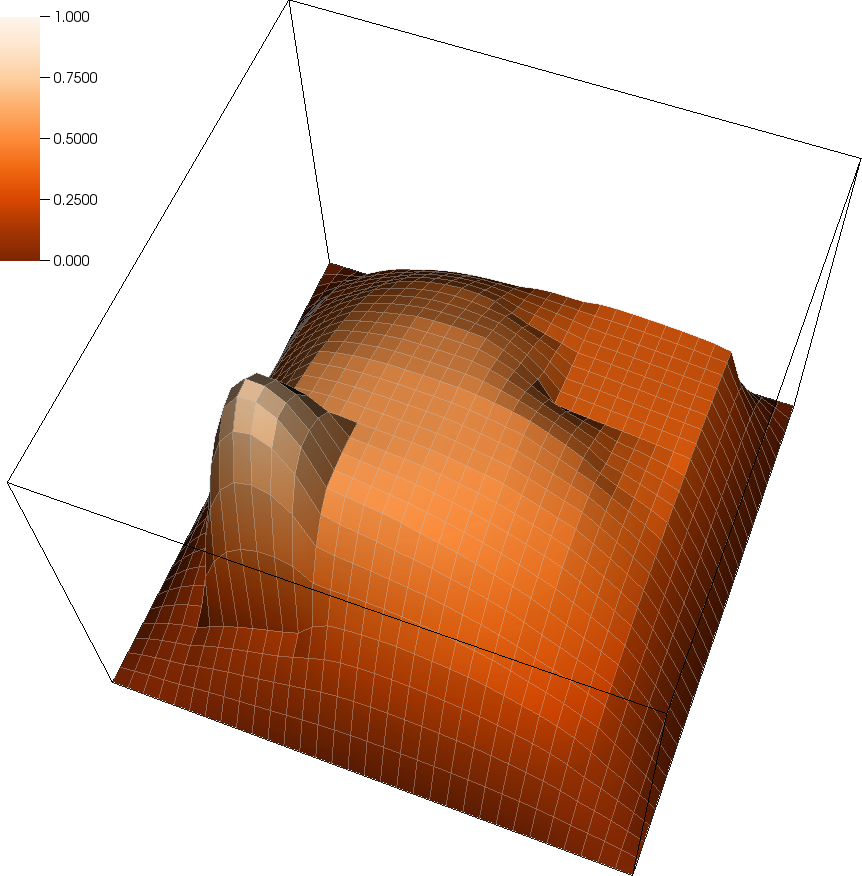}
  \end{tabular}
  \caption{\it Top left: The $8\times 8$ grid of values $\hat\theta$ used
    in generating the ``true measurements'' $\hat z$ via the forward
    model discussed in Section~\ref{sec:true-measurements}. The color
    scale matches that in Fig.~\ref{fig:mean}.
    Top right: The solution of the Poisson equation corresponding to
    $\hat\theta$ on the $256\times 256$ mesh and using $Q_3$ finite
    elements. ``True'' measurements
    $\hat z$ are obtained from this solution by evaluation at the 
    points shown in the right panel of Fig.~\ref{fig:forward}.
    Bottom left: For comparison, the solution obtained for the
    same set of parameters $\hat\theta$, but using the $32\times 32$ mesh and
    $Q_1$ element that defines the forward model. Bottom right:
    The solution of this discrete forward model applied to the posterior
    mean $\left<\theta\right>_{\pi(\theta|\hat z)}$ that we will compute later; the values of
    $\left<\theta\right>_{\pi(\theta|\hat z)}$
    are tabulated in Table~\ref{table:theta-means},
    and visualized in Fig.~\ref{fig:mean}.}
  \label{fig:hat-z}
\end{figure}

Indeed, this is how we have generated $\hat z$: We chose a set of
parameters $\hat\theta$ that corresponds to a membrane of uniform
stiffness $a(\mathbf x)=1$ except for two inclusions in which $a=0.1$
and $a=10$, respectively. This set up is shown in
Fig.~\ref{fig:hat-z}.%
\footnote{This set up has the accidental downside that both the
  set of parameters $\hat \theta$ and the set of measurement points
  $\mathbf x_k$ at which we evaluate the solution are symmetric
  about the diagonal of the domain. Since the same is true for
  our finite element meshes, the exact solution of the benchmark
  results in a probability distribution that is invariant
  to permutations of parameters about the diagonal as well, and
  this is apparent in Fig.~\ref{fig:mean}, for example. A better
  designed benchmark would have avoided this situation, but we only
  realized the issue after expending several years of CPU time. At the same time, the expected symmetry of values allows for a basic check of the correctness of inversion algorithms: If the inferred mean value $\left<\theta_7\right>_{\pi(\theta|\hat z)}$ is not approximately equal to $\left<\theta_{63}\right>_{\pi(\theta|\hat z)}$ -- see the numbering shown in the left panel of Fig.~\ref{fig:forward} -- then something is wrong.
  }
Using $\hat\theta$, we then used the series of
mappings as shown in \eqref{eq:forward-model} to compute $\hat
z$. However, to avoid an inverse crime, we have used a $256\times 256$
mesh and a bicubic ($Q_3$) finite element to compute 
$\hat\theta \mapsto \hat u_h \mapsto \hat z={\cal M}\hat u_h$, rather than the $32\times 32$ mesh and a bilinear
($Q_1$) element used to define the mapping $\theta\mapsto u_h \mapsto z^\theta={\cal M} u_h$.

As a consequence of this choice of higher accuracy (and higher
computational cost), we can in general not expect that there is a set
of parameters $\theta$ for which the forward model of
Section~\ref{sec:forward-model} would predict measurements $z^\theta$
that are \textit{equal} to $\hat z$. Furthermore, the presence of the prior
probability $\pi_\text{pr}$ in the definition of $\pi(\theta|\hat z)$
implies that we should not expect that either the mean
$\left<\theta\right>_{\pi(\theta|\hat z)}$ nor the MAP point
$\theta_\text{MAP}=\arg\max_\theta \pi(\theta|\hat z)$ are equal
or even just close to the ``true'' parameters $\hat\theta$.

\section{Statistical assessment of $\pi(\theta|\hat z)$}
\label{sec:assessment}
The previous section provides a concise definition of the
non-normalized posterior probability density $\pi(\theta|\hat
z)$. Given that the mapping $\theta\mapsto z^\theta$ is nonlinear and
involves solving a partial differential equation, there is no hope
that $\pi(\theta|\hat z)$ can be expressed as an explicit formula. On
the other hand, all statistical properties of $\pi(\theta|\hat z)$ can
of course be obtained by sampling, for example using algorithms such
as the Metropolis-Hastings sampler \cite{hastings1970monte}.

In order to provide a useful benchmark, it is necessary that at least
some properties of $\pi(\theta|\hat z)$ are known with sufficient
accuracy to allow others to compare the convergence of their sampling
algorithms. To this
end, we have used a simple Metropolis-Hastings sampler to compute
$\num{2e11}$ samples that characterize $\pi(\theta|\hat z)$, in the
form of $N=\num{2000}$ Markov chains of length $N_L=10^8$ each. (Details of
the sampling algorithm used to obtain these
samples are given in Appendix~\ref{sec:software-MH-sampler}.)
Using
the program discussed in Appendix~\ref{sec:software-forward}, the  effort to
produce this many samples amounts to approximately 30 CPU years on current
hardware. On the other hand, we will show below that this many samples
are really necessary in order to provide statistics of
$\pi(\theta|\hat z)$ accurately enough to serve as reference values -- at least, if one insists on
using an algorithm as simple as the Metropolis-Hastings method. In
practice, we hope that this benchmark is useful in the development of
algorithms that are substantially better than the Metropolis-Hastings
method. In addition, when assessing the convergence properties of a
sampling algorithm, it is of course not necessary to achieve the same
level of accuracy as we obtain here.

In the following let us therefore provide a variety of statistics
computed from our samples, along with an assessment of the accuracy
with which we believe that we can state these results. In the
following, we will denote by $0\le L<N=2000$ the number of the chain,
and $0\le \ell <N_L=10^8$ the number of a sample $\theta_{L,\ell}$ on
chain $L$. If we need to indicate one of the 64 components of a
sample, we will use a subscript index $k$ for this purpose as already used in Section~\ref{sec:forward-model}.

\subsection{How informative is our data set?}
\label{sec:information}

While we have $N=2000$ chains, each with a large number $N_L=10^8$ of samples per chain, 
a careful assessment needs to include an evaluation how informative all
of these samples really are. For example, if the samples on each chain
had a correlation length of $10^7$ because our Metropolis-Hastings
sampler converges only very slowly, then each chain really only contains
approximately ten statistically independent samples of $\pi(\theta|\hat z)$.
Consequently, we could not expect great accuracy in estimates of the
mean value, covariance matrices, and other quantities obtained from
each of the chains. Similarly, if the ``burn-in'' time of the sampler
is a substantial fraction of the chain lengths $N_L$, then we would
have to throw away many of the early samples.

To assess these questions, we have computed the autocovariance
matrices
\begin{align}
  AC_L(s) 
  &= 
  \left<[\theta_{L,\ell} - \left<\theta\right>_L]
        [\theta_{L,\ell-s} - \left<\theta\right>_L]^T
  \right>_L
  \notag
  \\
  \label{eq:autocovariance}
  &=
  \frac{1}{N_L-s-1}\sum_{\ell=s}^{N_L-1}
        [\theta_{L,\ell} - \left<\theta\right>_L]
        [\theta_{L,\ell-s} - \left<\theta\right>_L]^T
\end{align}
between samples $s$ apart on chain $L$. We expect samples with a small
lag $s$ to be highly correlated (i.e., $AC_L(s)$ to be a matrix that is large in some sense), whereas
for large lags $s$, samples should be uncorrelated and $AC_L(s)$ should
consequently be small. A rule of thumb 
is that samples at lags $s$ can be considered decorrelated from each other if $AC_L(s) \le 10^{-2} AC_L(0)$ entrywise; see Appendix C.

\begin{figure}
  \centering
  \includegraphics[width=0.46\textwidth]{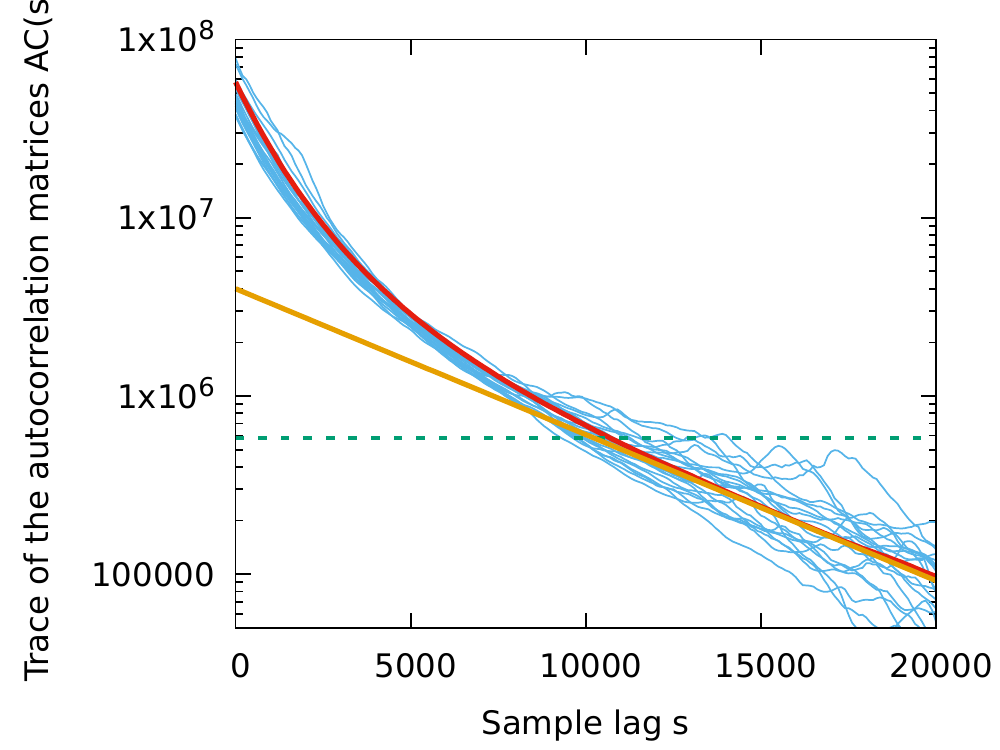}
  \caption{\it Decay of the trace of the autocovariance matrices $AC_L(s)$ with the
  sample lag $s$. The light blue curves show the traces of
  $AC_L(s)$ for twenty of our chains. The red curve shows the
  trace of the averaged autocovariance matrices, $AC(s)=\frac{1}{N}\sum_{L=0}^{N-1} AC_L(s)$.
  The thick yellow line corresponds to
  the function $4\cdot 10^6 \cdot e^{-s/5300}$ and shows the expected, asymptotically
  exponential decay of the autocovariance. The dashed green
  line indicates a reduction of the autocovariance by roughly
  a factor of 100 compared to its starting values.}
  \label{fig:autocorrelation}
\end{figure}

Fig.~\ref{fig:autocorrelation} shows the trace of these autocovariance matrices for several
of our chains. (We only computed the autocovariance at lags $s=0, 100, 200, 300,\ldots$ up to $s = 20,000$ because of the cost of computing $AC_L(s)$.)
The curves show that the
autocovariances computed from different chains all largely agree,
and at least asymptotically decay roughly exponentially with $s$ as expected.
The data also suggest that the autocorrelation length of our chains is
around $N_{AC}=10^4$ -- in other words, each of our chains should result
in approximately $N_L/N_{AC}=10^4$ meaningful and statistically 
independent samples.

To verify this claim, 
we estimated the {\em integrated autocovariance}~\cite{sokal1997monte} using
\begin{equation}\label{eq:IAT_matrix}
\IAC \approx   \frac{1}{N}\sum_{L=0}^{N-1}
100 \sum_{s=-200}^{200} AC_L(|100s|).
\end{equation}
The integrated autocovariance 
is obtained by summing up the autocovariance. (The factor of $100$ appears because we only computed $AC_L(s)$ at lags that are multiples of $100$.)  We show in Appendix C that the integrated autocovariance 
leads to the following
estimate of the effective sample size:
\begin{equation}\label{eq:ESS_perchain}
    \textup{Effective sample size per chain} \approx \frac{N_L}{\lambda_{\max}(C^{-1}\cdot \IAC)} \approx 1.3 \times 10^4,
\end{equation}
where $\lambda_{\max}$ indicates the 
maximum eigenvalue and $C = \frac{1}{N}\sum_{L=0}^{N-1}AC_L(0)$ is the covariance matrix; see also~\eqref{eq:covariance} below.
This is in good
agreement with Fig.~\ref{fig:autocorrelation} and the effective sample size derived from it.

\begin{figure}
  \includegraphics[width=0.46\textwidth]{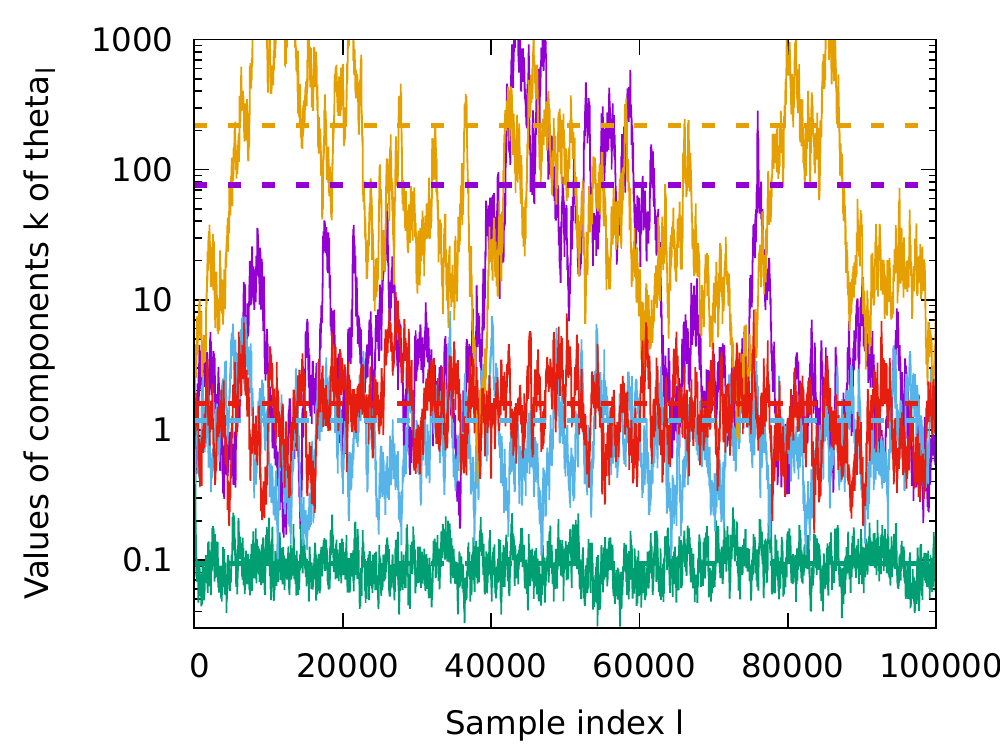}
  \hfill
  \includegraphics[width=0.46\textwidth]{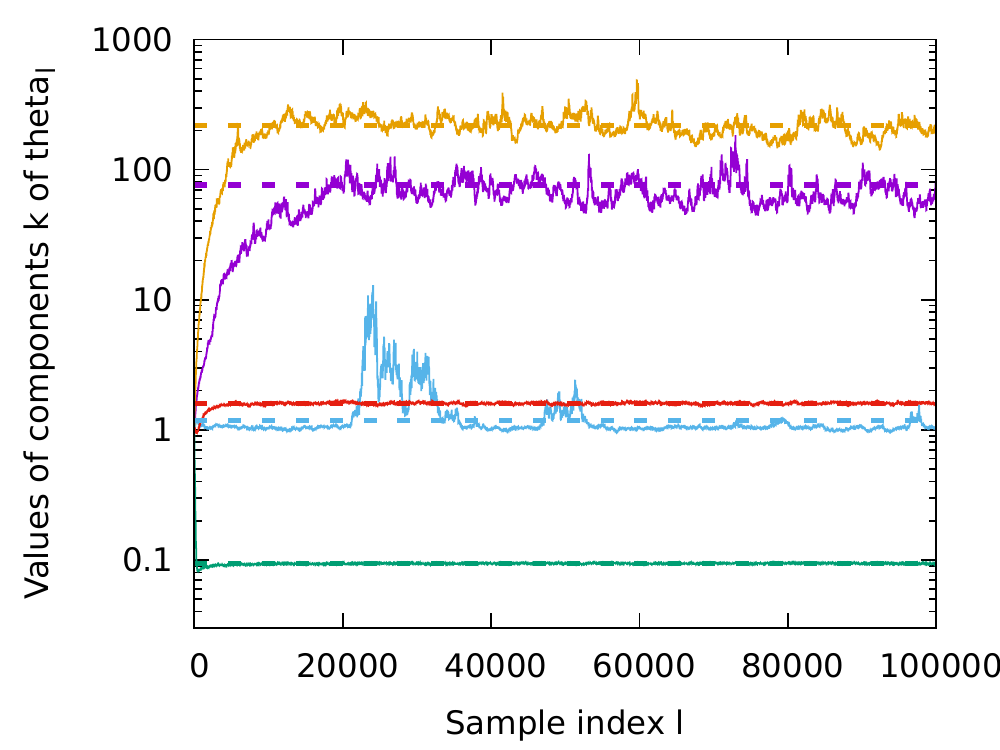}
  \caption{\it 
  Two perspectives on the length of the burn-in phase of our sampling scheme.
  Left: Values of components $k=0, 9, 36, 54$, and $63$ of 
  samples $\theta_\ell$ for $\ell=0,\ldots,100,000$, for one, randomly
  chosen chain. (For the geometric locations
  of the parameters that correspond to these components,
  see Fig.~\ref{fig:forward}.)
  Right: Across-chain averages
  $
    \frac{1}{N} \sum_{L=0}^{N-1} \theta_{L,\ell}
  $
  for the same components $k$ as above.
  Both images also shows the mean values
  $\left<\theta\right>_{\pi(\theta|\hat z)}$ for these five
  components as dashed lines. 
}
  \label{fig:burnin}
\end{figure}

This leaves the question of how long the burn-in period of our sampling scheme is. Fig.~\ref{fig:burnin} shows two perspectives on this question. The left panel of the figure shows several components $\theta_{\ell,k}$ of the samples
of one of our chains for the first few autocorrelation lengths. The data shows that there is at least no obvious ``burn-in'' period on this scale that would require us to throw away a substantial part of the chain. At the same time, it also illustrates that the large components of $\theta$ are poorly constrained and vary on rather long time scales that make it difficult to assess convergence to the mean. The right panel shows across-chain averages of the $\ell$th samples, more clearly illustrating that the burn-in period may only last around 20,000 samples -- that is, that only around 2 of the approximately 10,000 statistically independent samples of each chain are unreliable.

Having thus convinced ourselves that it is safe to use all $10^8$ samples
from all chains, and that there is indeed meaningful information
contained in them, we will next turn our attention towards computing
statistical information that characterizes $\pi(\theta|\hat z)$
and against which other implementations of sampling methods can
compare their results.

\subsection{The mean value of $\pi(\theta|\hat z)$}
\label{sec:mean}

\begin{figure}
  \centering
  \begin{tabular}{p{0.42\textwidth}p{0.42\textwidth}}
    \includegraphics[width=0.35\textwidth]{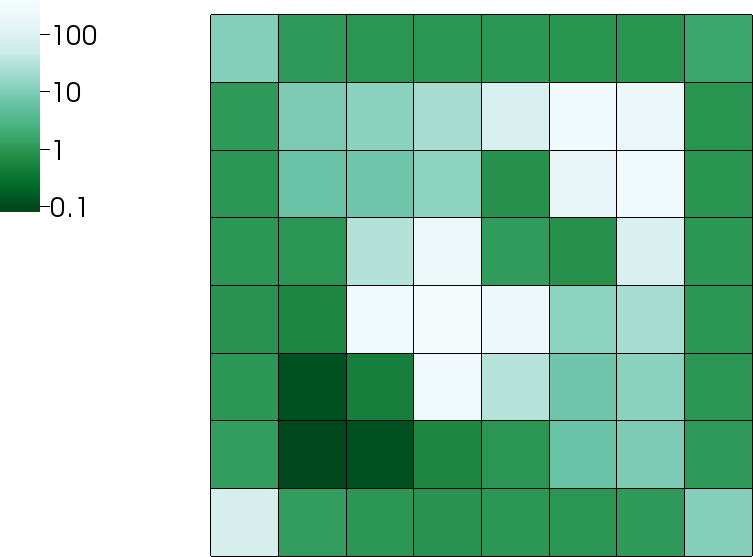}
    &
    \includegraphics[width=0.35\textwidth]{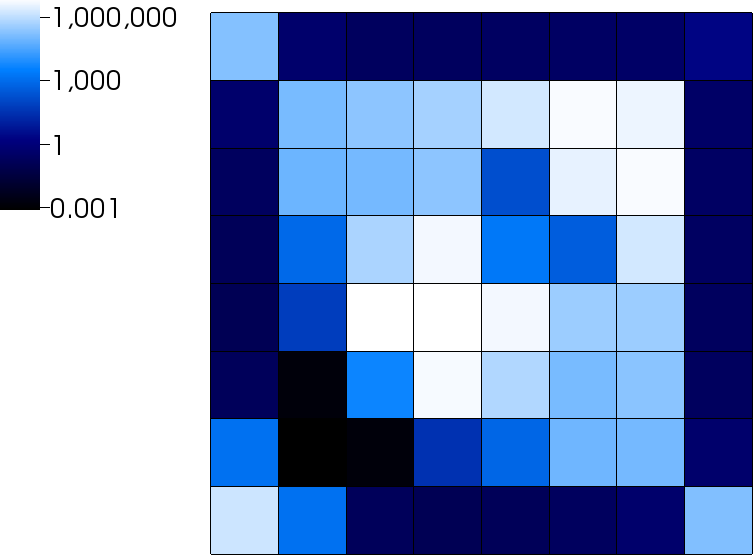}
    \\
    \centering \small Inferred mean values $\left<\theta_k\right>_{\pi(\theta|\hat z)}$.
    &
    \centering \small Variances $C_{kk}$.
  \end{tabular}
  \caption{\it Left: Visualization on an $8\times 8$ grid of mean values
  $\left<\theta\right>_{\pi(\theta|\hat z)}$ obtained from
  our samples. This figure should be compared with the values
  $\hat\theta$ used as input to generate the ``true measurements''
  $\hat z$ and shown in Fig.~\ref{fig:hat-z}. 
  Right: Variances $C_{kk}$ of
  the parameters. Dark colors indicate that a parameter is accurately known; light colors that the variance is large.
  Standard deviations (the square roots of the variances) are larger
  than the mean values in some cases because of the
  heavy tails in the distributions of
  parameters.}
  \label{fig:mean}
\end{figure}

The simplest statistic one can compute from samples of a distribution
is the mean value. Table~\ref{table:theta-means} shows the 64 values that
characterize the mean
\begin{align*}
  \left<\theta_k\right>_{\pi(\theta|\hat z)}
  =
  \frac{1}{N}\sum_{L=0}^{N-1}\left(\frac{1}{N_L}\sum_{\ell=0}^{N_L-1}
  \theta_{L,\ell,k}\right).
\end{align*}
A graphical
representation of $\left<\theta\right>_{\pi(\theta|\hat z)}$ is shown Fig.~\ref{fig:mean}
and can be compared to the ``true'' values $\hat\theta$ shown in
Fig.~\ref{fig:hat-z}.%
\footnote{By comparing Fig.~\ref{fig:hat-z} and the data of
  Table~\ref{table:theta-means} and Fig.~\ref{fig:mean}, it is clear that for some parameters, the mean $\left<\theta_k\right>_{\pi(\theta|\hat z)}$ is far away from the value $\hat\theta_k$ used to generate the original data $\hat z$ -- principally for those parameters that correspond to large values $\hat\theta_k$, but also the ``white cross'' below and left of center. For the first of these two places, we can first note that the prior probability $\pi_\text{pr}$ defined in \eqref{eq:prior} is quite heavy-tailed, with a mean far larger than where its maximum is located. And second, by realizing that a membrane that is locally very stiff is not going to deform in a substantially different way in response to a force from one that is even stiffer in that region -- in other words, in areas where the coefficient $a^\theta(\mathbf x)$ is large, the likelihood function \eqref{eq:likelihood} is quite insensitive to the exact values of $\theta$, and the posterior probability will be dominated by the prior $\pi_\text{pr}$ with its large mean.

For the ``white cross'', one can make plausible that the likelihood is
uninformative and that, consequently, mean value and variances are
again determined by the prior. To understand why this is so, one could
imagine by analogy what would happen if one could measure the solution
$u(\mathbf x)$ of \eqref{eq:laplace-1}--\eqref{eq:laplace-2} exactly
and everywhere, instead of only at a discrete set of points. In that
case, we would have $u(\mathbf x)=z(\mathbf x)$, and we could infer
the coefficient $a(\mathbf x)$ by solving
\eqref{eq:laplace-1}--\eqref{eq:laplace-2} \textit{for the
  coefficient} instead of for $u$. This leads to the advection-reaction equation $-\nabla
z(\mathbf x) \cdot \nabla a(\mathbf x) - (\Delta z(\mathbf x))
a(\mathbf x) = f(\mathbf x)$, which is ill-posed and does not provide
for a stable solution $a(\mathbf x)$ at those places where $\nabla
z(\mathbf x)=\nabla u(\mathbf x)\approx 0$. By comparison with
Fig.~\ref{fig:hat-z}, we can see that at the location of the white
cross, we could not identify the coefficient at one point even if we
had measurements available \textit{everywhere}, and not stably so in
the vicinity of that point. We can expect that this is also so in the
discrete setting of this benchmark -- and that consequently, at this
location, only the prior provides information.}

\begin{table}
  \caption{\it Sample means $\left<\theta\right>_{\pi(\theta|\hat z)}$ for the 64 parameters,
    along with their estimated 2-sigma uncertainties.}
   
  \tiny
  \centering
  \begin{minipage}[t]{0.33\textwidth}
  \begin{tabular}{l@{ }r@{ $\pm$ }l}
$\left<\theta\right>_{0}=$ 
             &\num{76.32}         &\num{0.30}\\
                &\num{1.2104}       &\num{0.0094}\\
         &\num{0.977380}        &\num{0.000051} \\
         &\num{0.882007}       &\num{0.000039} \\
         &\num{0.971859}       &\num{0.000048} \\
            &\num{0.947832}      &\num{0.000064} \\
            &\num{1.08529}       &\num{0.00011}\\
             &\num{11.39}           &\num{0.10}\\
$\left<\theta\right>_{8}=$ 
              &\num{1.119}        &\num{0.011}\\
        &\num{0.0937215}      &\num{0.0000027} \\
             &\num{0.1157992}      &\num{0.0000039}\\
              &\num{0.5815}       &\num{0.0022}\\
         &\num{0.9472}       &\num{0.0079}\\
          &\num{6.258}        &\num{0.079}\\
              &\num{9.334}        &\num{0.090}\\
              &\num{1.08151}      &\num{0.00011}\\
$\left<\theta\right>_{16}=$
        &\num{0.977449}      &\num{0.000052}\\
              &\num{0.1157962}      &\num{0.0000038}\\
              &\num{0.461}        &\num{0.020}\\
             &\num{267.01}           &\num{0.55}\\
              &\num{30.87}         &\num{0.19}\\
          &\num{7.189}        &\num{0.089}\\
              &\num{12.39}          &\num{0.11}\\
         &\num{0.949863}      &\num{0.000073}\\
  \end{tabular}
  \end{minipage}
  \begin{minipage}[t]{0.32\textwidth}
  \begin{tabular}{l@{}r@{ $\pm$ }l}
$\left<\theta\right>_{24}=$
        &\num{0.881977}      
        &\num{0.000039} \\
         &\num{0.5828}       &\num{0.0020}\\
              &\num{267.72}         &\num{0.62}\\
                &\num{369.35}         &\num{0.64}\\
                 &\num{234.59}          &\num{0.53}\\
              &\num{13.29}         &\num{0.14}\\
          &\num{22.36}         &\num{0.16}\\
              &\num{0.988806}      &\num{0.000074}\\
$\left<\theta\right>_{32}=$
        &\num{0.971900}       &\num{0.000049}\\
              &\num{0.9509}       &\num{0.0079}\\
              &\num{30.76}         &\num{0.19}\\
              &\num{233.93}         &\num{0.52}\\
                &\num{1.169}        &\num{0.012}\\
         &\num{0.8327}       &\num{0.0057}\\
          &\num{88.52}          &\num{0.33}\\
         &\num{0.987809}      &\num{0.000079}\\
$\left<\theta\right>_{40}=$
             &\num{0.947816}      &\num{0.000065}\\
                &\num{6.260}        &\num{0.076}\\
          &\num{7.119}        &\num{0.087}\\
              &\num{13.20}         &\num{0.13}\\
              &\num{0.8327}       &\num{0.0035}\\
              &\num{176.73}         &\num{0.44}\\
              &\num{283.38}         &\num{0.58}\\
         &\num{0.914212}      &\num{0.000077}\\
  \end{tabular}
  \end{minipage}
  \begin{minipage}[t]{0.32\textwidth}
  \begin{tabular}{l@{}r@{ $\pm$ }l}
$\left<\theta\right>_{48}=$
              &\num{1.08521}      &\num{0.00011}\\
          &\num{9.386}       &\num{0.089}\\
              &\num{12.44}         &\num{0.12}\\
                &\num{22.50}         &\num{0.17}\\
          &\num{88.57}         &\num{0.33}\\
              &\num{283.41}         &\num{0.57}\\
              &\num{218.65}         &\num{0.49}\\
                &\num{0.933451}      &\num{0.000087}\\
$\left<\theta\right>_{56}=$
              &\num{11.35}         &\num{0.11}\\
              &\num{1.08143}      &\num{0.00011}\\
              &\num{0.949869}      &\num{0.000074}\\
         &\num{0.988770}      &\num{0.000074}\\
         &\num{0.987866}      &\num{0.000083}\\
              &\num{0.914247}      &\num{0.000077}\\
         &\num{0.933426}      &\num{0.000087}\\
              &\num{1.59984}     &\num{0.00030}\\
                          \\
                          \\
                          \\
                          \\
                          \\
                          \\
                          \\
                          \\
  \end{tabular}  
  \end{minipage}
  \label{table:theta-means}
\end{table}

To assess how accurately we know this
average, we consider that we have $N=2000$ chains of length
$N_L=10^8$ each, and that each of these has its own chain averages
\begin{align*}
  \left<\theta_k\right>_L
  =
  \frac{1}{N_L}\sum_{\ell=0}^{N_L-1}
  \theta_{L,\ell,k}.
\end{align*}
The ensemble average $\left<\theta_k\right>_{\pi(\theta|\hat z)}$ is of
course the average of the chain averages $\left<\theta_k\right>_L$ across chains, but
the chain averages vary between themselves and we can compute the
standard deviation of these chain averages as
\begin{align*}
  \text{stddev}\left(\left<\theta_k\right>_L\right)
  =
  \left[
  \frac{1}{N}
  \sum_{L=0}^{N-1} \left(\left<\theta_k\right>_L -
  \left<\theta_k\right>_{\pi(\theta|\hat z)} \right)^2
  \right]^{1/2}.
\end{align*}
Under standard assumptions, and assuming that the posterior is Gaussian,
we can then estimate that we know the
ensemble averages $\left<\theta_k\right>_{\pi(\theta|\hat z)}$ to
within an accuracy of $\pm \frac{1}{\sqrt{N}}
\text{stddev}\left(\left<\theta_k\right>_L\right)$
with 68\% (1-sigma) certainty, and with an accuracy
of $\pm \frac{2}{\sqrt{N}}
\text{stddev}\left(\left<\theta_k\right>_L\right)$
with 95\% (2-sigma) certainty. In reality, the posterior is \textit{not} Gaussian (see Section~\ref{sec:nongaussianty}), and the argument is not true as stated; however, computing 2-sigma values for all parameters is still a useful metric for how accurately we know each of the parameters.

This 2-sigma accuracy is also
provided in  Table~\ref{table:theta-means}. For all but parameter $\theta_{18}$ (for which the relative 2-sigma uncertainty is $4.5\%$), the relative uncertainty in $\left<\theta\right>$ is between $0.003\%$ and $1.3\%$. 
In other words, the table provides nearly two certain digits for all but one parameter, and four digits for at least half of all parameters.

\subsection{The covariance matrix of $\pi(\theta|\hat z)$ and its properties}
\label{sec:covariance}
The second statistic we demonstrate is the covariance matrix,
\begin{align}\begin{split}\label{eq:covariance}
  C_L
  &=
  \frac{1}{N_L-1}
  \sum_{\ell=0}^{N_L-1}
  \left(\theta_{L,\ell}-\left<\theta\right>_{\pi(\theta|\hat z)}\right)
  \left(\theta_{L,\ell}-\left<\theta\right>_{\pi(\theta|\hat z)}\right)^T,
  \\
  C &= \frac{1}{N}\sum_{L=0}^{N-1} C_L.
  \end{split}
\end{align}
While conceptually easy to compute, in practice it is substantially harder to obtain accuracy in $C$ than it is to compute accurate means $\left<\theta\right>_{\pi(\theta|\hat z)}$: While we know the latter to two or more digits of accuracy, see Table~\ref{table:theta-means}, there is substantial variation between the matrices $C_L$.%
\footnote{For diagonal entries $C_{L,kk}$, the standard deviation of the variation \textit{between} chains is between 0.0024 and 37.7 times the corresponding entry $C_{kk}$ of the average covariance matrix. The variation can be even larger for the many small off-diagonal entries. On the other hand, the average (across chains) difference $\|C_L-C\|_F$ is $0.9\|C\|_F$. This would suggest that we don't know very much about these matrices, but as shown in the rest of the section, \textit{qualitative} measures can be extracted robustly.}
The remainder of this section therefore only provides qualitative conclusions we can draw from our estimate of the covariance matrix, rather than providing quantitive numbers.

\begin{figure}
  \centering
  \includegraphics[width=0.44\textwidth]{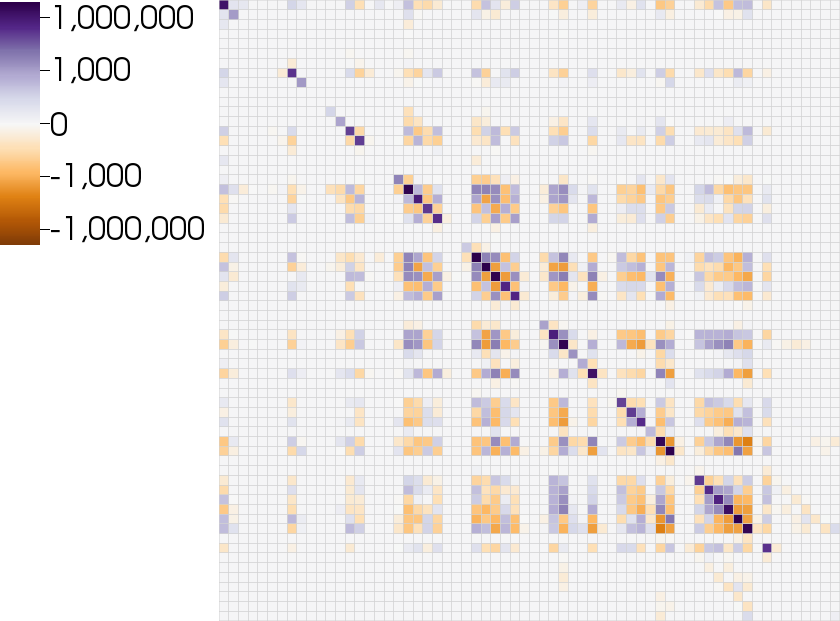}
  \hfill
  \includegraphics[width=0.44\textwidth]{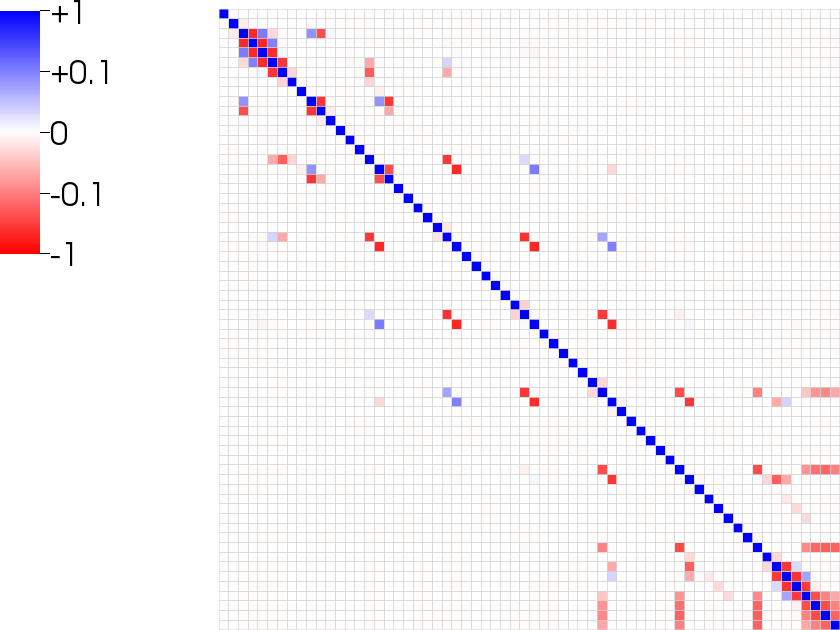}
  
  \caption{\it Left: A visualization of the covariance matrix $C$
    computed from the posterior probability distribution $\pi(\theta|\hat z)$. The sub-structure of the matrix in the form of
    $8\times 8$ tiles represents that geographically neighboring -- and consequently correlated -- parameters
    are either a distance of $\pm 1$ or $\pm 8$ apart. Right: Correlation matrix 
    $D_{ij}=\frac{C_{ij}}{\sqrt{C_{ii}}\sqrt{C_{jj}}}$.}
  \label{fig:covariance}
\end{figure}

First, the diagonal entries of $C$, $C_{kk}$, provide the variances of the
statistical distribution of $\theta_k$, and are shown on the right of Fig.~\ref{fig:mean}; the off-diagonal
entries $C_{k\ell}$ suggest how \textit{correlated} parameters $\theta_k$ and
$\theta_l$ are and are depicted in Fig.~\ref{fig:covariance}.

In the context of inverse problems related to partial differential
equations, it is well understood that we expect the parameters to
be highly correlated. This can be understood intuitively given that we
are thinking of a membrane model: If we increased the stiffness value
on one of the $8\times 8$ pixels somewhat, but decreased the stiffness
value on a neighboring correspondingly, then we would expect to obtain
more or less the same global deformation pattern -- maybe there are
small changes at measurement points close to the perturbation, but for
measurement points far away the local perturbation will make little
difference. As a consequence, we 
should expect that $L(\hat z|\theta)\approx L(\hat z|\tilde\theta)$
where $\theta$ and $\tilde\theta$ differ only in two nearby components, one
component of $\tilde\theta$ being slightly larger and the other being
slightly smaller than the corresponding component of $\theta$. If the
changes are small, then we will also have that $\pi(\theta|\hat
z)\approx \pi(\tilde\theta|\hat z)$ -- in other words, we would expect
that $\pi$ is approximately constant in the \textit{secondary diagonal
  directions} $(0,\ldots,0,+\varepsilon,0,\ldots,0,-\varepsilon,0,\ldots
0)$ in $\theta$ space.

On the other hand, increasing (or decreasing) the stiffness value in
\textit{both} of two adjacent pixels just makes the membrane overall
more (or less) stiff, and will yield different displacements at
\textit{all} measurement locations. Consequently, we expect that the
posterior probability distribution $\pi(\theta|\hat z)$ will
strongly vary in the \textit{principal diagonal
  directions} $(0,\ldots,0,+\varepsilon,0,\ldots,0,+\varepsilon,0,\ldots
  0)$ in $\theta$ space.

\begin{figure}
  \centering
  \phantom{.}
  \hfill
  \includegraphics[width=0.44\textwidth]{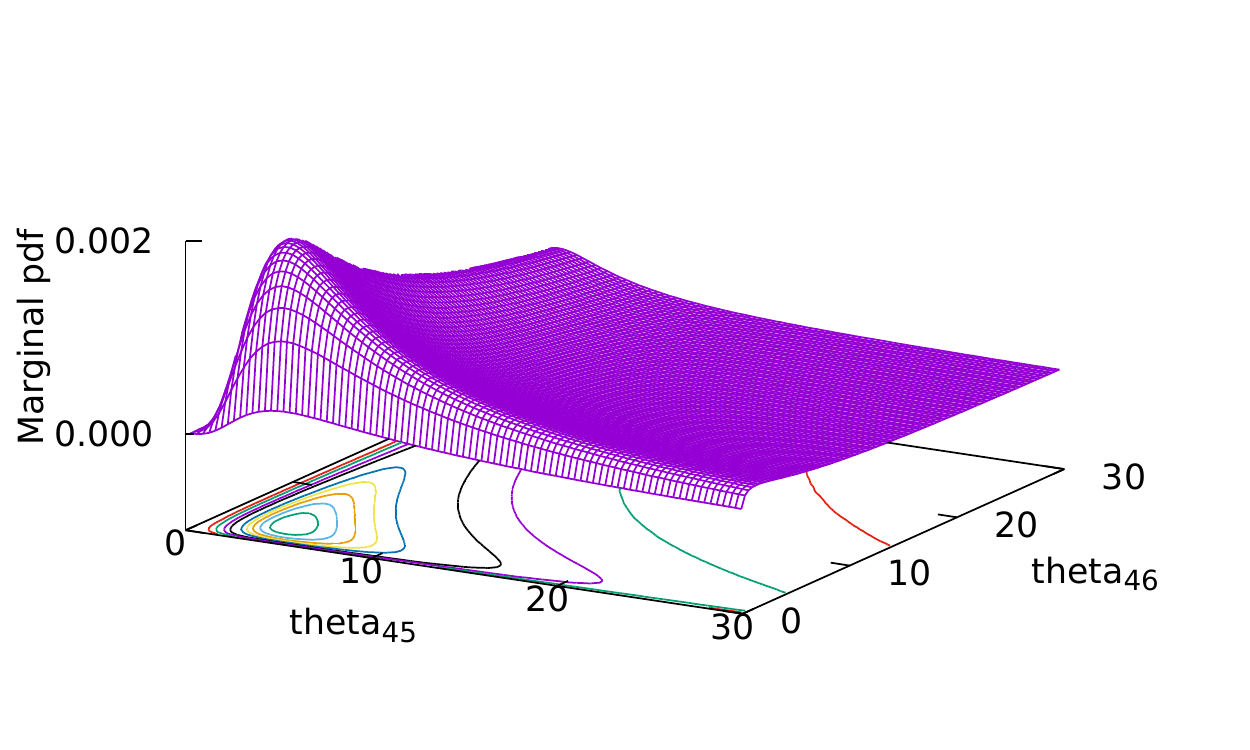}
  \hfill
  \includegraphics[width=0.44\textwidth]{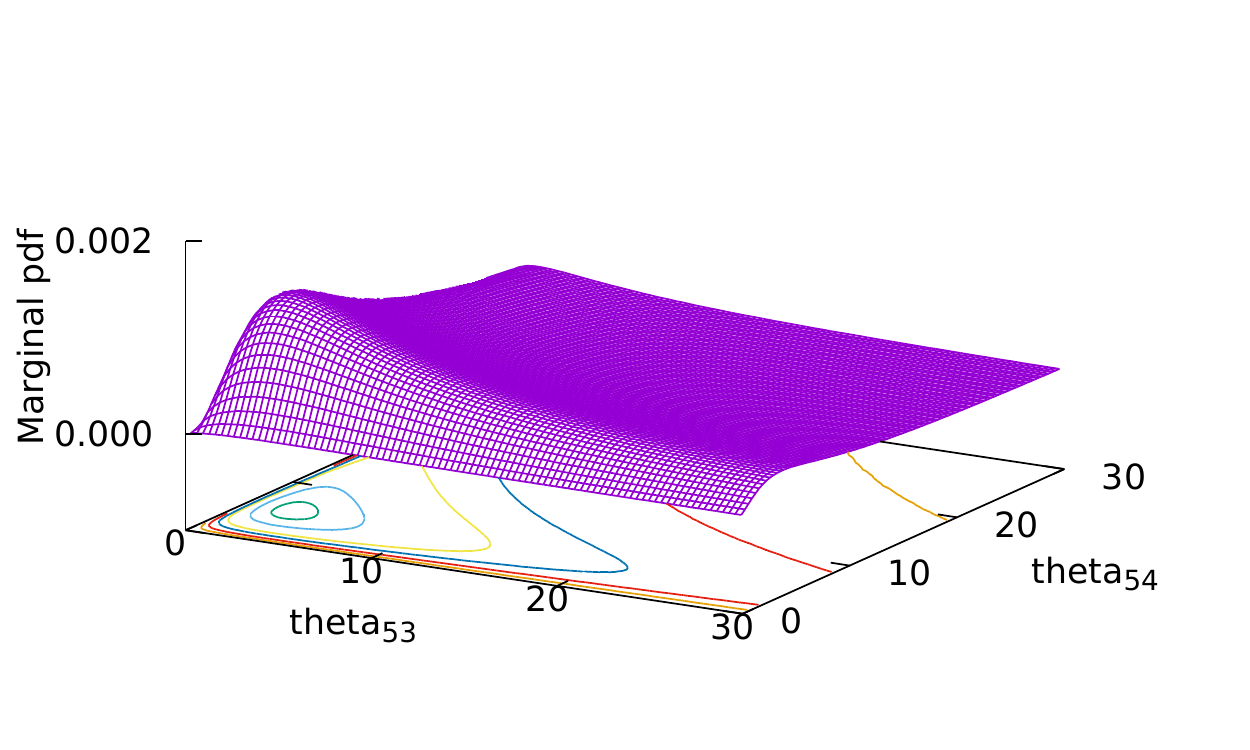}
  \hfill
  \phantom{.}

  \caption{\it Pairwise marginal probabilities for $\theta_{45}$ and $\theta_{46}$ (left) and for 
  $\theta_{53}$ and $\theta_{54}$ (right), respectively. See Fig.~\ref{fig:forward} for the relative locations of these parameters.
  These marginal distributions illustrate the anti-correlation of parameters: If one is large, the other is most likely small, and vice versa.}
  \label{fig:pair-histogram}
\end{figure}

We can illustrate this by computing two-dimensional histograms of the
samples for parameters $\theta_k$ and $\theta_l$ corresponding to
neighboring pixels -- equivalent to a two-dimensional marginal
distribution. We show such histograms in Fig.~\ref{fig:pair-histogram}. These also
indicate that the posterior probability distribution $\pi(\theta|\hat
z)$ is definitely not Gaussian -- see also Remark~\ref{remark:sigma}.

\begin{figure}
  \centering
  \includegraphics[width=0.8\textwidth]{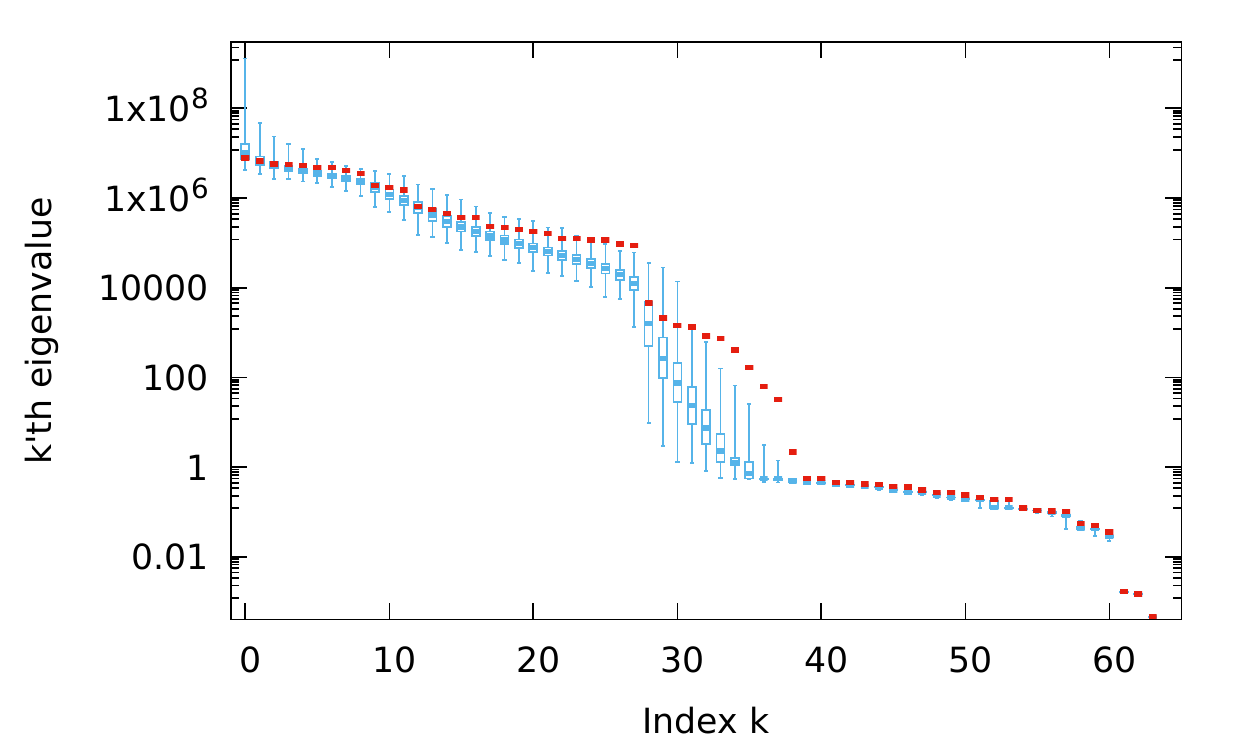}

  \caption{\it Singular values of the covariance matrix $C$ as defined in \eqref{eq:covariance}. Red bars indicate eigenvalues of the across-chain averaged covariance matrix. Blue boxes correspond to the 25th and 75th percentiles of the corresponding eigenvalues of the covariance matrices $C_L$ of individual chains; vertical bars extend to the minimum and maximum across chains for the $k$th eigenvalue of the matrices $C_L$; blue bars in the middle of boxes indicate the median of these eigenvalues.}
  \label{fig:singular-values}
\end{figure}

A better way to illustrate correlation is to compute a singular value
decomposition of the covariance matrix $C$. Many inverse problems have only a
relatively small number of large singular values of $C$
\cite{hippylib,Flath2011,BuiThanh2012,Worthen2014,Petra2014,chen2019}, 
suggesting that only
a finite number of modes is resolvable with the data
available -- in other words, the problem is ill-posed. Fig.~\ref{fig:singular-values} shows the singular values of the covariance
matrix $C$ for the current case. The data suggests that from the 169 measured pieces of (noisy) data, a deterministic inverse problem could only recover some 25-30 modes of the parameter vector $\theta\in \R^{64}$ with reasonable accuracy.%
\footnote{The figure shows the spread of each of the eigenvalues of the within-chain matrices $C_L$ in blue, and the eigenvalues of the across-chain matrix $C$ in red. One would expect the latter to be well approximated by the former, and that is true for the largest and smallest eigenvalues, but not for the ones in the middle. There are two reasons for this: First, each of the $C_L$ is nearly singular, but because each chain is finite, the poorly explored directions are different from one chain to the next. At the same time, it is clear that the sum of (different) singular matrices may actually be ``less singular'', with fewer small eigenvalues, and this is reflected in the graph.  A second reason is that we computed the eigenvalues of each of the $C_L$ and ordered them by size when creating the plot, but without taking into account the associated eigenspaces. As a consequence, if one considers the, say, 32nd largest eigenvalue of $C$, the figure compares it with the 32nd largest eigenvalues of all of the $C_L$, when the correct comparison would have been with those eigenvalues of the matrices $C_L$ whose eigenspace is most closely aligned; this may be an eigenvalue elsewhere in the order, and the effect will likely be the most pronounced for those eigenvalues whose sizes are the least well constrained.

The conclusions to be drawn from Fig.~\ref{fig:singular-values} are therefore not the actual sizes of eigenvalues, but the number of ``large'' eigenvalues. This observation is robust, despite the inaccuracies in our determination of $C$.
}

\subsection{Higher moments of $\pi(\theta|\hat z)$}
\label{sec:nongaussianty}

In some sense, solving Bayesian inverse problems is not very interesting
if the posterior distribution $p(\theta|\hat z)$ for the parameters is Gaussian,
or at least approximately so, because it can be done much more efficiently by computing the maximum likelihood estimator through a deterministic inverse problem, and then computing the covariance matrix via the Hessian of the deterministic (constrained) optimization problem.
For example, \cite{hippylib} provides an excellent overview of the techniques that can be used in this case.
Because of these simplifications, it is of interest to know how close the posterior density of this benchmark is to a multi-dimensional Gaussian.

\begin{figure}
  \centering
  \includegraphics[width=0.6\textwidth]{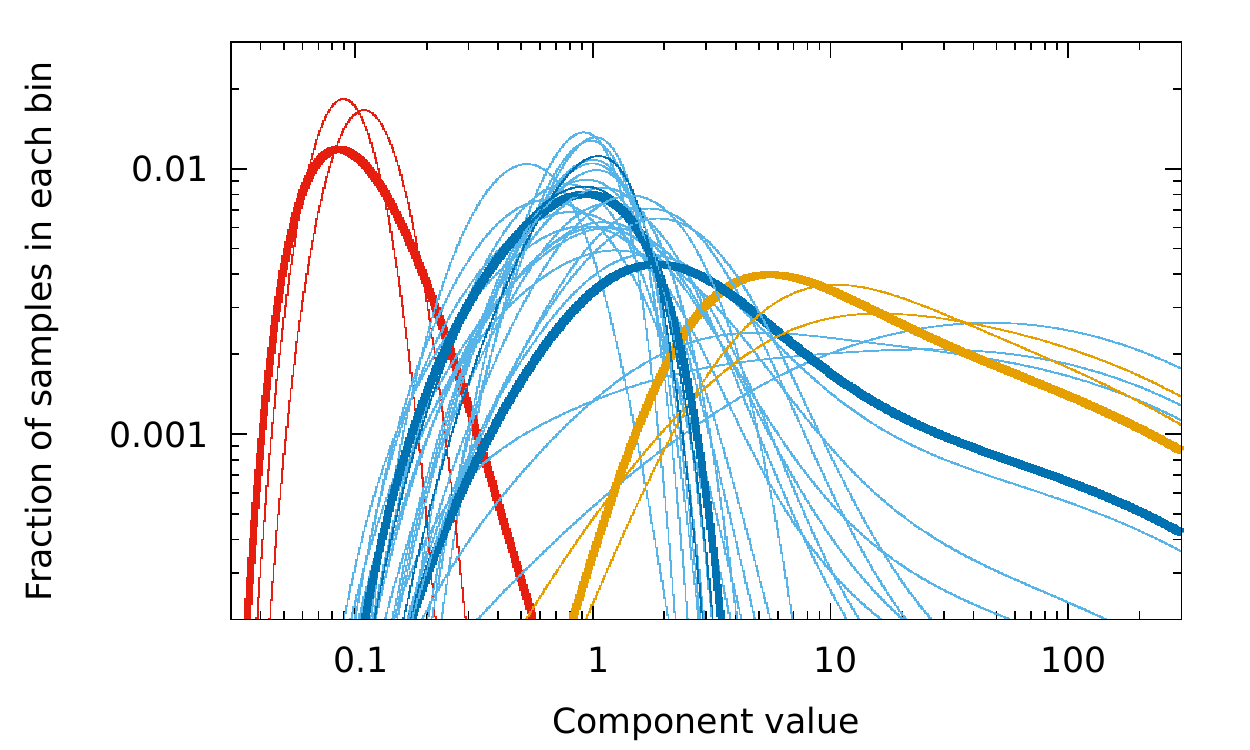}

  \caption{\it Histograms (marginal distributions) for all 64 components of $\theta$, 
    accumulated over all \num{2e11} samples. The
    histograms
    for those parameters whose values were 0.1 for the purposes
    of generating the vector $\hat z$ are shown in red, those
    whose values were 10 in orange, and all others in blue.
    The figure highlights histograms for some of those components
    $\theta_k$ whose marginal distributions are clearly neither
    Gaussian nor log-Gaussian.
    Note that the histograms were generated with bins
    whose size increases exponentially from left to
    right, so that they are depicted with equal size in
    the figure given the logarithmic $\theta_k$ axis. 
    Histograms with equal bin size on a linear scale
    would be more heavily biased towards the left, but with
    very long tails to the right. See Fig.~\ref{fig:pair-histogram} for (pair) histograms using a linear scale.}
  \label{fig:histograms}
\end{figure}

To evaluate this question, Fig.~\ref{fig:histograms} shows histograms
of all of the parameters, using 1000 bins that are equally spaced 
in logarithmic space; i.e., for each component $k$, we create 1000 bins between -3 and +3 and
sort samples into these bins based on $\log_{10}(\theta_k)$. It is clear that many of the parameters have heavy
tails and can, consequently, not be approximated well by Gaussians.
On the other hand, given the prior distribution \eqref{eq:prior}
we attached to each
of the parameters, it would make sense to conjecture that the
\textit{logarithms} $\log(\theta_k)$ might be Gaussian distributed.

If that were so, the double-logarithmic plot shown in the figure
would consist of histograms in the form of parabolas open to the bottom,
and again a simpler -- and presumably cheaper to compute --
representation of $\pi(\theta|\hat z)$ would be possible. However,
as the figure shows, this too is clearly not the case: While some
parameters seem to be well described by such a parabola, many others
have decidedly non-symmetric histograms, or shapes that are simply not
parabolic. As a consequence, we conclude that the benchmark is at
least not boring in the sense that its posterior distribution could
be computed in some comparably much cheaper way.

\subsection{Rate of convergence to the mean}
\label{sec:rate-of-convergence}

The data provided in the previous subsections allows checking whether 
a separate implementation of this benchmark converges to the same
probability distribution $\pi(\theta|\hat z)$. However, it does not
help in assessing whether it does so faster or slower than the
simplistic Metropolis-Hastings method used herein. Indeed,
as we will outline in the next section, we hope that this work
spurs the development and evaluation of methods that can achieve
the same results \textit{without} needing more than $10^{11}$ samples.

To this end, let us here provide metrics for how fast
our method converges to the mean 
$\left<\theta\right>_{\pi(\theta|\hat z)}$ discussed in Section~\ref{sec:mean}.
More specifically, if we denote by
$\left<\theta\right>_{L,n}=\frac{1}{n}\sum_{\ell=0}^{n-1} \theta_{L,\ell}$ the running mean of samples zero to $n-1$ on
chain $L$, then we are interested in how fast it converges to
the mean. We measure this using the following error norm
\begin{align}
    e_L(n)
    &=
    \left\|
    \text{diag}\,(\left<\theta\right>_{\pi(\theta|\hat z)})^{-1}
    \left( \left<\theta\right>_{L,n}
          -\left<\theta\right>_{\pi(\theta|\hat z)}\right)
    \right\|
    \notag
    \\
    &=
    \left[
    \left( \left<\theta\right>_{L,n}
          -\left<\theta\right>_{\pi(\theta|\hat z)}\right)^T
    \text{diag}\,(\left<\theta\right>_{\pi(\theta|\hat z)})^{-2}
    \left( \left<\theta\right>_{L,n}
          -\left<\theta\right>_{\pi(\theta|\hat z)}\right)
    \right]^{1/2}.
    \label{eq:error-norm}
\end{align}
The weighting by a diagonal matrix containing the inverses of the
estimated parameters $\left<\theta\right>_{\pi(\theta|\hat z)}$ (given in Table~\ref{table:theta-means} and known to sufficient accuracy)
ensures that
the large parameters with their large variances do not
dominate the value of $e_L(n)$. In other words, $e_L(n)$ corresponds to the ``root mean squared relative error''.%
\footnote{A possibly better choice for the weighting would be to use the 
  inverses of the diagonal entries of the covariance matrix -- i.e., the
  variances of the recovered marginal probability distributions
  of each parameter. However, these are only approximately known -- see
  the discussion in Section~\ref{sec:covariance} -- and consequently
  do not lend themselves for a concise definition of a reproducible
  benchmark.}

\begin{figure}
  \centering
  \includegraphics[width=0.6\textwidth]{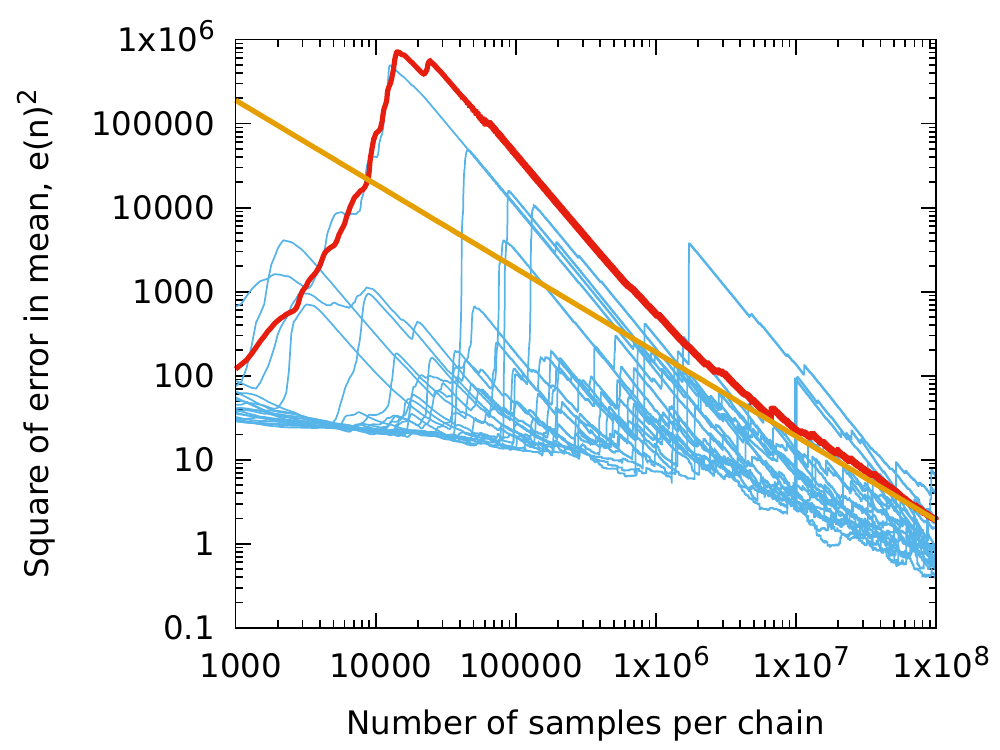}
  \caption{\it Convergence of the square of the relative error $e_L(n)^2$ between the running mean
    up to sample $n$ and the true mean
    $\left<\theta\right>_{\pi(\theta|\hat z)}$, measured in the
    weighted norm \eqref{eq:error-norm}, for a subset
    of chains. The thick red curve corresponds to the average of these squared errors over all chains. This average squared error is dominated by some chains with large errors and so lies above the curves $e_L(n)^2$ of most chains.
    The thick orange line corresponds to the decay $\frac{\num{1.9e8}}{n}$ and represents an approximately average convergence behavior of the chains.}
  \label{fig:convergence}
\end{figure}

Fig.~\ref{fig:convergence} shows the convergence of a few chains to
the ensemble average. While there is substantial variability between chains, it is clear that for each chain $e_L(n)^2\rightarrow 0$ and, furthermore, that this convergence follows the classic
one-over-$n$ convergence of statistical sampling algorithms.
Indeed, averaging $e_L(n)^2$ over all chains,
\begin{align*}
    e(n)^2
    =
    \frac{1}{N} \sum_{L=0}^{N-1} e_L(n)^2,
\end{align*}
the behavior of this decay of the ``average'' square error $e(n)^2$ can be approximated by the
following formula that corresponds to the orange line in the figure:
\begin{align}
  \label{eq:convergence}
    e(n)^2 \approx \frac{\num{1.9e8}}{n}.
\end{align}
While we have arrived at the factor $\num{1.9e8}$ by fitting a curve ``by eye'', it turns out -- maybe remarkably -- that
we can also theoretically support this behavior: using the Markov chain central limit theorem~\cite{jones2004markov} (see Appendix C for details), 
we can estimate the mean of $ne(n)^2$ by
\begin{equation*}
\textup{tr}\left(\diag(\langle\theta\rangle_{\pi(\theta|\hat{z})})^{-1}\cdot \IAC \cdot
\diag(\langle\theta\rangle_{\pi(\theta|\hat{z})})^{-1}\right)\approx 1.9 \times 10^8,
\end{equation*}
where the matrix $\IAC$ is defined in~\eqref{eq:IAT_matrix}.

If one measures computational effort by how many times an algorithm
evaluates the probability distribution $\pi(\theta|\hat z)$, then
$n$ in \eqref{eq:error-norm} can be interpreted as work units
and \eqref{eq:convergence} provides an approximate relationship
between work and error. Similar relationships can be obtained
experimentally for other sampling algorithms that we hope this
benchmark will be used for, and 
\eqref{eq:convergence} therefore allows comparison
among other algorithms as well as against the one used here.

\section{Conclusions and what we hope this benchmark achieves}
\label{sec:hope}

As the data presented in the previous section illustrates, it is
possible to obtain reasonably accurate statistics about the Bayesian solution
of the benchmark introduced in Section~\ref{sec:benchmark-description},
even using a rather simple method: The standard Metropolis-Hastings
sampler. At the same time, using this method, it is not at all trivial to
compute posterior statistics \textit{accurately}: We had to compute
$\num{2e11}$ samples, and expended 30 CPU years on this task (plus another
two CPU years on postprocessing the samples).

But all of this also makes this a good benchmark: Simple algorithms, with
known performance, \textit{can} solve it to a reasonable accuracy, and more
advanced algorithms should be able to do so in a fraction of time without
making the test case trivial. For example, it is not unreasonable to hope
that advanced sampling software \cite{Debusschere2017,Dakota,MUQ,McDougall2017}, using
multi-level and multi-fidelity expansions \cite{Peherstorfer2018,Dodwell2015,Seo2020,Fleeter2020}, and maybe in conjunction with methods that exploit the structure of the problem to approximate covariance matrices \cite{hippylib}, might be able to reduce the compute
time by a factor of 100 to 1000, possibly also running computations in
parallel. This would move characterizing the performance of such
algorithms for the case at hand to the range of a few hours or days on moderately parallel computers; practical computations might not actually need the same level of
accuracy and could be solved even more rapidly.

As a consequence of these considerations, we hope that providing a benchmark
that is neither too simple nor too hard, and for which the solution is known
to good accuracy, spurs research in the development of better sampling
algorithms for Bayesian inverse problems. Many such algorithms of course
already exist, but in many cases, their performance is not characterized
on standardized test cases that would allow a fair comparison. In particular,
their performance is often characterized using probability distributions
whose characteristics have nothing to do with those that result from
inverse problems -- say, sums of Gaussians. By
providing a standardized benchmark that matches what we expect to see in
actual inverse problems -- along with an open source
implementation of a code that computes the posterior probability function
$\pi(\theta|\hat z)$ (see Appendix~\ref{sec:software}) -- we hope that
we can contribute to more informed comparisons between newly proposed
algorithms: Specifically, that their performance can be compared with
the relationship shown in \eqref{eq:convergence} and Fig.~\ref{fig:convergence} to
provide a concrete factor of speed-up over the method used here.

\subsection*{Acknowledgments}

W.~Bangerth was partially
supported by the National Science Foundation under award OAC-1835673
as part of the Cyberinfrastructure for Sustained Scientific Innovation (CSSI)
program; by award DMS-1821210; 
by award EAR-1925595;
and by the Computational Infrastructure
in Geodynamics initiative (CIG), through the National Science
Foundation under Award No.~EAR-1550901 and The
University of California -- Davis.

W.~Bangerth also gratefully acknowledges the discussions and early
experiments with Kainan Wang many years ago. These early attempts
directly led to the ideas that were ultimately encoded in this
benchmark. He also appreciates the collaboration with 
Mantautas Rimkus and Dawson Eliasen on the \textsc{SampleFlow} library that was used
for the statistical evaluation of samples. Finally, the feedback given
by Noemi Petra, Umberto Villa, and Danny Long are acknowledged with gratitude.

D. Aristoff gratefully acknowledges support from the National Science Foundation via awards DMS-1818726 
and DMS-2111277.

\appendix

\section{An open source code to sample $\pi(\theta|\hat z)$}
\label{sec:software}
We make a code that implements this benchmark available as part of the
``code gallery'' for \dealii{} at {\small \url{https://dealii.org/developer/doxygen/deal.II/code_gallery_MCMC_Laplace.html}}, using the name
\texttt{MCMC-Laplace}, and using the Lesser GNU Public License (LGPL)
version 2.1 or later as the license.
\dealii{} is a software library that provides the basic tools
and building blocks for writing finite element codes that solve
partial differential equations numerically. More information about
\dealii{} is available at \cite{arndt2019dealii,dealII92}. The \dealii{} code
gallery is a collection of programs based on \dealii{} that were
contributed by users as starting points for others' experiments.

The code in question has essentially three parts: (i) The forward
solver that, given a set of parameters $\theta$, produces the output
$z^\theta$ using the map discussed in Section~\ref{sec:forward-model};
(ii) the statistical model that implements the likelihood
$L(z|\theta)$ and the prior probability $\pi_\text{pr}(\theta)$, and
combines these to the posterior probability $\pi(\theta|\hat z)$; and
(iii) a simple Metropolis-Hastings sampler that draws samples from
$\pi(\theta|\hat z)$. The second of these pieces is quite trivial,
encompassing only a couple of functions; we will therefore only
comment on the first and the third piece below.

\subsection{Details of the forward solver}
\label{sec:software-forward}
The forward solver is a C++ class whose main function performs the
following steps:
\begin{enumerate}
\item It takes a 64-dimensional vector $\theta$ of parameter values,
  and interprets it as the coefficients that describe a piecewise
  constant field $a(\mathbf x)$;
\item It assembles a linear system that corresponds to the finite
  element discretization of equations
  \eqref{eq:laplace-1}--\eqref{eq:laplace-2} using a $Q_1$ (bilinear)
  element on a uniformly refined mesh;
\item It solves this linear system to obtain the solution vector
  $U^\theta$ that corresponds to the function $u_h^\theta(\mathbf x)$; and
\item It evaluates the solution $u_h^\theta$ at the measurement points
  $\mathbf x_k$ to obtain $z^\theta$.
\end{enumerate}
It then returns $z^\theta$ to the caller for evaluation with the
statistical model.

Such a code could be written in \dealii{} with barely more than 100
lines of C++ code, and this would have been sufficient for the purpose
of evaluating new ideas of sampling methods. However, we wanted to
draw as large a number of samples as possible, and consequently
decided to see how fast we can make this code.

To this end, we focused on accelerating three of the operations listed
above, resulting in a code that can evaluate $\pi(\hat z|\theta)$ in
4.5 \si{ms} on an Intel Xeon E5-2698 processor with 2.20GHz (on which
about half of the samples used in this publication were computed),
3.1 \si{ms} on an AMD EPYC 7552 processor with 2.2 \si{GHz} (the other half),
and 2.7 \si{ms} on an Intel Core i7-8850H
processor with 2.6 \si{GHz} in one of the authors' laptops.

The first part of the code that can be optimized for the current
application uses the fact that the linear system that needs to be
assembled is the sum of contributions from each of the cells of the
mesh. Specifically, the contribution from cell $K$ is
\begin{align*}
  A^K = P_K^T A^K_\text{local} P_K
\end{align*}
where $P_K$ is the restriction from the global set of degrees of
freedom to only those degrees of freedom that live on cell $K$, and
$A_\text{local}^K$ is -- for the $Q_1$ element used here -- a $4\times
4$ matrix of the form
\begin{align*}
  (A^K_\text{local})_{ij}
  =
  \int_K a^\theta(\mathbf x)
  \nabla \varphi_i(\mathbf x) \cdot \nabla \varphi_j(\mathbf x)
  \, \text{d}x.
\end{align*}
This suggests that the assembly, including the integration above that
is performed via quadrature, has to be repeated every time we
consider a new set of parameters $\theta$. However, since we
discretize $u_h$ on a mesh that is a strict refinement of the one used
for the coefficient $a^\theta(\mathbf x)$, and because $a^\theta(\mathbf x)$
is piecewise constant, we can note that
\begin{align*}
  (A^K_\text{local})_{ij}
  =
  \theta_{k(K)}
  \underbrace{
    \int_K 
    \nabla \varphi_i(\mathbf x) \cdot \nabla \varphi_j(\mathbf x)
    \, \text{d}x}_{(A_\text{local})_{ij}},
\end{align*}
where $k(K)$ is the index of the element of $\theta$ that corresponds
to cell $K$. Here, the matrix $A_\text{local}$ no longer depends on
$\theta$ and can, consequently, be computed once and for all at the
beginning of the program. Moreover, $A_\text{local}$
does not actually depend on the cell $K$ as long as all cells have the
same shape, as is the case here. We therefore have to store only one
such matrix. This approach makes the assembly substantially faster
since we only have to perform the local-to-global operations
corresponding to $P_K$ on every cell for every new $\theta$, but no
longer any expensive integration/quadrature.

Secondly, we have experimented with solving the linear systems so
assembled as fast as possible. For the forward solver used for each
sample, the size of these linear systems is $1089\times 1089$, with at
most 9 entries per row. Following a substantial amount of experimentation, we found that a sparse direct
solver is faster than any of the other approaches we have tried, and
we use the UMFPACK \cite{umfpack} interfaces in \dealii{} for this
purpose. In particular, this approach is faster than attempting to use
an algebraic multigrid method as a solver or preconditioner for the
Conjugate Gradient method. We have also tried to use a sparse
decomposition via UMFPACK as a preconditioner for the CG method,
updating the decomposition only every few samples -- based on the
assumption that the samples change only relatively slowly and so a
decomposition of the matrix for one sample is a good preconditioner
for the matrix corresponding to a subsequent sample. However, this
turned out to be slower than using a new decomposition for each
sample.

The linear solver described above consumes about 90\% of the time necessary to
evaluate each sample. As a consequence, there is certainly further
room for improvements. After all numerical results had been generated for this publication, we have followed up on a suggestion by Martin Kronbichler to replace the linear solver by a conjugate gradient method preconditioned by an incomplete LU decomposition. This accelerates the computations by about a factor of three, from \SI{2.7}{ms} to less than \SI{0.9}{ms} per sample on the Intel Core i7-8850H processor mentioned above. It is this accelerated version that is available at the website mentioned above; we have verified that the new version results in the same results up to at least 11 digits in computing the posterior probability using the techniques mentioned in Appendix~\ref{sec:testing}.

Finally, evaluating the solution of a finite element field
$u_h(\mathbf x)$ at arbitrary points $\mathbf x_k$ is an expensive
operation since one has to find which cell $K$ the point belongs to
and then transform this point into the reference coordinate system of
the cell $K$. On the other hand, the point evaluation is a linear and bounded
operation, and so there must exist a vector $m_k$ so that $u_h(\mathbf
x_k)=m_k \cdot U$ where $U$ is the vector of coefficients that
describe $u_h$. This vector $m_k=(\varphi_i(\mathbf
x_k))_{i=0}^{1089-1}$ can be computed once and for all. The
computation of $z=(u_h(\mathbf x_k))_{k=0}^{169-1}$ can then be
facilitated by building a matrix $M$ whose rows are the vectors
$m_k$, and then the evaluation at all measurement points reduces to
the operation $z=MU$. $M$ is a sparse matrix with at most 4 entries
per row, making this a very economical approach.

The code with all of these optimizations is not very large --
it contains 197 semicolons.%
\footnote{Counting semicolons is a commonly used metric in C and C++
  programs. It roughly coincides with the number of declarations and
  statements in a program, and is a better metric for \textit{code
    size} than the number of lines of code, as the latter also
  includes comments and empty lines used to help the readability of a
  code.}

\subsection{Details of the Metropolis-Hastings sampler}
\label{sec:software-MH-sampler}

The steps described at the start of this appendix yield an algorithm that, given a sample $\theta$,
can evaluate $\pi(\theta|\hat z)$. We use this functionality to drive
a Metropolis-Hastings sampler to obtain a large number of samples
characterizing the posterior probability distribution.

While the basic algorithm of the Metropolis-Hastings sampler is well
known \cite{hastings1970monte,Kaipio2005}, its practical implementation depends crucially on a
number of details that we will describe in the following.

First, we start the sampling process with a fixed sample
$\theta_0=(1,\ldots,1)^T$, corresponding to a coefficient
$a^\theta(\mathbf x)=1$.

Secondly, an important step in the Metropolis-Hastings algorithm is
the generation of a ``trial'' sample $\tilde{\theta}_{\ell}$ based on the
current sample $\theta_{\ell}$. To this end, we use the following strategy: We define the components of $\tilde{\theta}_{\ell,k}$ of $\tilde{\theta}_{\ell}$ as
\begin{align*}
  \tilde{\theta}_{\ell,k}= e^{\ln(\theta_{\ell,k}) + \xi_k}
    = \theta_{\ell,k} e^{\xi_k},
\end{align*}
where $\xi_k$, $k=0,\ldots,63$, are iid Gaussians with mean 
$0$ and standard deviation 
$\sigma_\text{prop}$.
In other 
words, the ``proposal distribution'' for the trial samples is an isotropic
Gaussian ball centered at $\theta_\ell$ in log space. This has the effect
that all elements of samples always stay positive, as one would
expect given that they correspond to material stiffness
coefficients. The use of a ball in log space is also consistent with
the description of our prior probability distribution in
Section~\ref{sec:prior}, which is also defined in log space.

To compute the 
Metropolis-Hastings acceptance probability, we first need 
to compute the proposal 
probability density. By definition
\begin{align*}
{\mathbb P}(\tilde{\theta}_{\ell,k} \le \tilde{t}_k | \theta_{\ell,k} = t_k) &=
{\mathbb P}(t_ke^{\xi_k} \le \tilde{t}_k) \\
&= {\mathbb P}(\xi_k \le \log \tilde{t}_k-\log t_k) \\
&= 
\frac{1}{\sqrt{2\pi\sigma_{\textup{prop}}^2}}\int_{-\infty}^{\log \tilde{t}_k-\log t_k} \exp\left(\frac{-x^2}{2\sigma_{\textup{prop}}^2}\right)\,dx.
\end{align*}
The probability density $p_{\textup{prop}}(\tilde{\theta}_{\ell,k}|\theta_{\ell,k})$ 
of proposing $\tilde{\theta}_{\ell,k}$ 
given $\theta_{\ell,k}$ is then the derivative of this expression with respect to $\tilde{t}_k$, with $\tilde{\theta}_{\ell,k}$ and ${\theta}_{\ell,k}$ in place of $\tilde{t}_k$ and $t_k$:
\begin{align*}
p_{\textup{prop}}(\tilde{\theta}_{\ell,k}|\theta_{\ell,k}) = \frac{1}{\tilde{\theta}_{\ell,k}}\frac{1}{\sqrt{2\pi\sigma_{\textup{prop}}^2}}\exp\left(\frac{-(\log \tilde{\theta}_{\ell,k} - \log \theta_{\ell,k})^2}{2\sigma_{\textup{prop}}^2}\right).
\end{align*}
By definition, the 
components $\tilde{\theta}_{\ell,k}$ of the proposal vector $\tilde{\theta}_{\ell}$ are 
independent conditional 
on the current state $\theta_{\ell}$. 
Thus the joint probability density, 
$p_{\textup{prop}}(\tilde{\theta}_{\ell}|\theta_{\ell})$, 
of proposing vector $\tilde{\theta}_{\ell}$ 
given vector $\theta_{\ell}$ is the 
product of the probabilities above:
\begin{equation*}
p_{\textup{prop}}(\tilde{\theta}_{\ell}|\theta_{\ell}) = \frac{1}{\prod_{k=0}^{63} \tilde{\theta}_{\ell,k}} \times\frac{1}{(2\pi\sigma_{\textup{prop}}^2)^{32}}\exp\left(-\frac{\sum_{k=0}^{63} (\log \tilde{\theta}_{\ell,k} - \log \theta_{\ell,k})^2}{2\sigma_{\textup{prop}}^2}\right).
\end{equation*}
The Metropolis-Hastings 
acceptance probability, $A(\tilde{\theta}_{\ell}|\theta_{\ell})$, to accept proposal 
$\tilde{\theta}_{\ell}$ given 
the current state $\theta_{\ell}$, is then
\begin{align*}
A(\tilde{\theta}_{\ell}|\theta_{\ell}) &=  \min\left\{1,\frac{\pi(\tilde{\theta}_{\ell}|\hat{z})}{\pi({\theta}_{\ell}|\hat{z})}\times\frac{p_{\textup{prop}}({\theta}_{\ell}|\tilde{\theta}_{\ell})}{p_{\textup{prop}}(\tilde{\theta}_{\ell}|\theta_{\ell})}\right\} \\
&= \min\left\{1,\frac{\pi(\tilde{\theta}_{\ell}|\hat{z})}{\pi({\theta}_{\ell}|\hat{z})}\times 
\prod_{k=0}^{63} \frac{\tilde{\theta}_{\ell,k}}{{\theta}_{\ell,k}}\right\}.
\end{align*}
As usual, with probability $A(\tilde{\theta}_{\ell}|\theta_{\ell})$ the 
proposal $\tilde{\theta}_{\ell}$ 
is accepted, in which 
case it 
becomes the next sample $\theta_{\ell+1} = \tilde{\theta}_{\ell}$; 
otherwise the proposal is 
rejected and we keep the current 
sample, 
$\theta_{\ell+1} = {\theta}_{\ell}$.

In our experiments, we use $\sigma_\text{prop}=0.0725$, corresponding
to changing $\theta$ by a factor $e^\xi$ that with 65\% probability
lies within the range $[e^{-\sigma_\text{prop}},e^{+\sigma_\text{prop}}]=[0.93,1.075]$. This results
in an acceptance rate for the Metropolis-Hastings algorithm of just under
24\%. This is close to the recommended value of 0.234 for 
Metropolis-Hastings sampling algorithms that can be
derived for specific probability distributions that
are generally simpler than the one we are interested in
here
\cite{Gelman97,Sherlock2009};
the theory guiding the derivation of the 0.234 value
may not be applicable here (see \cite{Bdard2008}),
but absent better guidance, we stuck with it.

Finally, all steps in the Metropolis-Hastings algorithms that require
the use of a random number use the MT19937 random number generator, as
implemented by C++11 compilers in the \texttt{std::mt19937} class.

\subsection{Implementations of the benchmark in alternative languages}

We consider the C++ implementation discussed above as the
``reference implementation'' of the benchmark. However, we recognize that it is rather heavy-weight in the sense that it requires the installation of the \textsc{deal.II} library. While this readily facilitates otherwise non-trivial modifications (e.g., for multilevel sampling schemes that require the solution of the forward problem on coarser meshes), it is clear that for simple experiments, it would be nice to have stand-alone implementations of the benchmark.

As a consequence, we have also developed Matlab and Python versions of the benchmark. These are available from the same website from which the C++ implementation is available -- see the link at the top of the appendix. These alternative implementations provide everything one needs to build a sampler: Namely, the functionality to provide an input vector $\theta$ and to compute the prior $\pi_\text{pr}(\theta)$ and the likelihood $L(\hat z|\theta)$ from such an input. Using \eqref{eq:posterior}, these can then be used to compute the posterior probability $\pi(\theta|\hat z)$ that forms the basis of most sampling algorithms. The Matlab 
version includes a basic 
Metropolis-Hastings sampler 
with parallel functionality.

On a recent laptop, the Matlab version is able to compute the posterior probability for a sample in about \SI{4}{ms}, not substantially slower than the C++ version we have used for our results. The Python version -- although an almost literal transcription of the Matlab version -- requires approximately \SI{25}{ms}. One imagines that it could be optimized further (for example, a substantial part of the run time is spent in the insertion of the cell-local $4\times 4$ matrices into the global matrix), but we have not attempted to do so.

\subsection{Testing of alternative implementations}
\label{sec:testing}

In order to facilitate testing of alternative implementations such as the ones discussed in the previous section (or testing modifications made to the C++ implementation itself), the website from which the benchmark can be obtained (see the top of this appendix) also contains a set of known input/output pairs. Specifically, it contains files for ten different input vectors $\theta$, along with the corresponding outputs $z^\theta$, likelihood $L(\hat z| \theta)$, and prior $\pi_\text{pr}(\theta)$ for each input vector. The latter two can be combined via \eqref{eq:posterior} to the posterior probability $\pi(\theta|\hat z)$ associated with the input vector $\theta$.

We have used these known input/output pairs obtained from our reference implementation to verify that the alternative implementations of the benchmark discussed in the previous sub-section are correct. For example, the Matlab implementation provides the vectors $z^\theta$ to relative errors on the order of $10^{-13}$; log priors and likelihoods are computed to relative errors less than $10^{-11}$. The Python version achieves the same level of accuracy.

\section{One dimensional version of the benchmark}
\label{sec:1d}
Many of the features of the posterior probability on the
64-dimensional parameter space that we have experimentally
observed in Section~\ref{sec:assessment} match those that
one would expect for inverse problems of the kind
discussed herein. In particular, the fact that we have a large
spread between the large and small eigenvalues of the
covariance matrix, the non-Gaussianity of the probability
distribution, and the anticorrelation of parameters defined on
neighboring cells did not come as a surprise. Yet, strict
proofs for these properties are hard to come by.

At the same time, we can investigate an analogous situation
in a one-dimensional model with parameters $\theta_k$, for which one can derive the
posterior probability distribution analytically.
In this Appendix, we work out some details 
of a one-dimensional version 
of the benchmark. Consider the 
generalized Poisson equation 
\begin{align}
-\frac{d}{dx}\left(a(x)\frac{du}{dx}(x)\right) &= f(x), &\qquad \qquad \qquad &0<x<1, \label{eq:poiss1D_1}\\
u(x) &= 0, &\qquad \qquad &x \in \{0,1\}.\label{eq:poiss1D_2}
\end{align} 
Again we assume $a(x) = a^\theta(x)$ 
is parametrized by
$\theta = (\theta_0,\ldots,\theta_{N-1})$ 
where 
\begin{equation}\label{eq:a_1D}
a(x) = \theta_k, \qquad \text{if}\quad \frac{k}{N} < x < \frac{k+1}{N},
\end{equation}
and we take $f(x) \equiv 1$.
The solution to~\eqref{eq:poiss1D_1} is then of the form
\begin{equation}\label{eq:poisson1D_sol}
u(x) = - \int_0^x \frac{y+c}{a(y)}\,dy, 
\end{equation}
where $c$ is constant determined by requiring $u(1)=0$.
Due to the piecewise constant form of $a = a(x)$, this rewrites as 
\begin{equation}\label{eq:piecewise}
u(x) = u_k(x) := -\frac{1}{\theta_k}\left(\frac{x^2}{2} + cx\right) + d_k, \qquad \text{if} \quad \frac{k}{N}<x<\frac{k+1}{N}.
\end{equation}
There are $N-1$ continuity constraints and two boundary conditions, namely, 
\begin{align}\label{eq:1D_constraints}
&u_{k-1}\left(\frac{k}{N}\right) = u_{k}\left(\frac{k}{N}\right), \qquad k=1,\ldots,N-1,\\
&u_0(0) = 0, \quad u_{N-1}(1) = 0.
\end{align}
This translates, via~\eqref{eq:piecewise}, into $N+1$ linear equations for the $N+1$ coefficients
$c$,$d_0,\ldots,d_{N-1}$.
The boundary conditions show that 
\begin{equation}
d_0 = 0, \qquad c = -\frac{\int_0^1 \frac{y}{a(y)}\,dy }{\int_0^1 \frac{1}{a(y)}\,dy} = -\frac{1}{2N}\frac{\sum_{k=0}^{N-1}\theta_k^{-1}\left(2k+1\right)}{\sum_{k=0}^{N-1} \theta_k^{-1}}. 
\end{equation}
The remaining $N-1$ equations for $d_1,\ldots,d_{N-1}$ can be solved using the continuity constraints, which give the equations
\begin{equation}
d_{k-1}-d_k = \left(\frac{1}{2}\left(\frac{k}{N}\right)^2+ c\frac{k}{N}\right)\left(\theta_{k-1}^{-1}-\theta_k^{-1}\right), \qquad k=1,\ldots,N-1,
\end{equation}
which have solution 
\begin{equation*}
d_k = -\sum_{j=1}^k \left(\frac{1}{2}\left(\frac{j}{N}\right)^2+ c\frac{j}{N}\right)\left(\theta_{j-1}^{-1}-\theta_j^{-1}\right), \qquad k=1,\ldots,N-1.
\end{equation*}

\begin{figure}
\centering
\includegraphics[width=0.6\textwidth]
{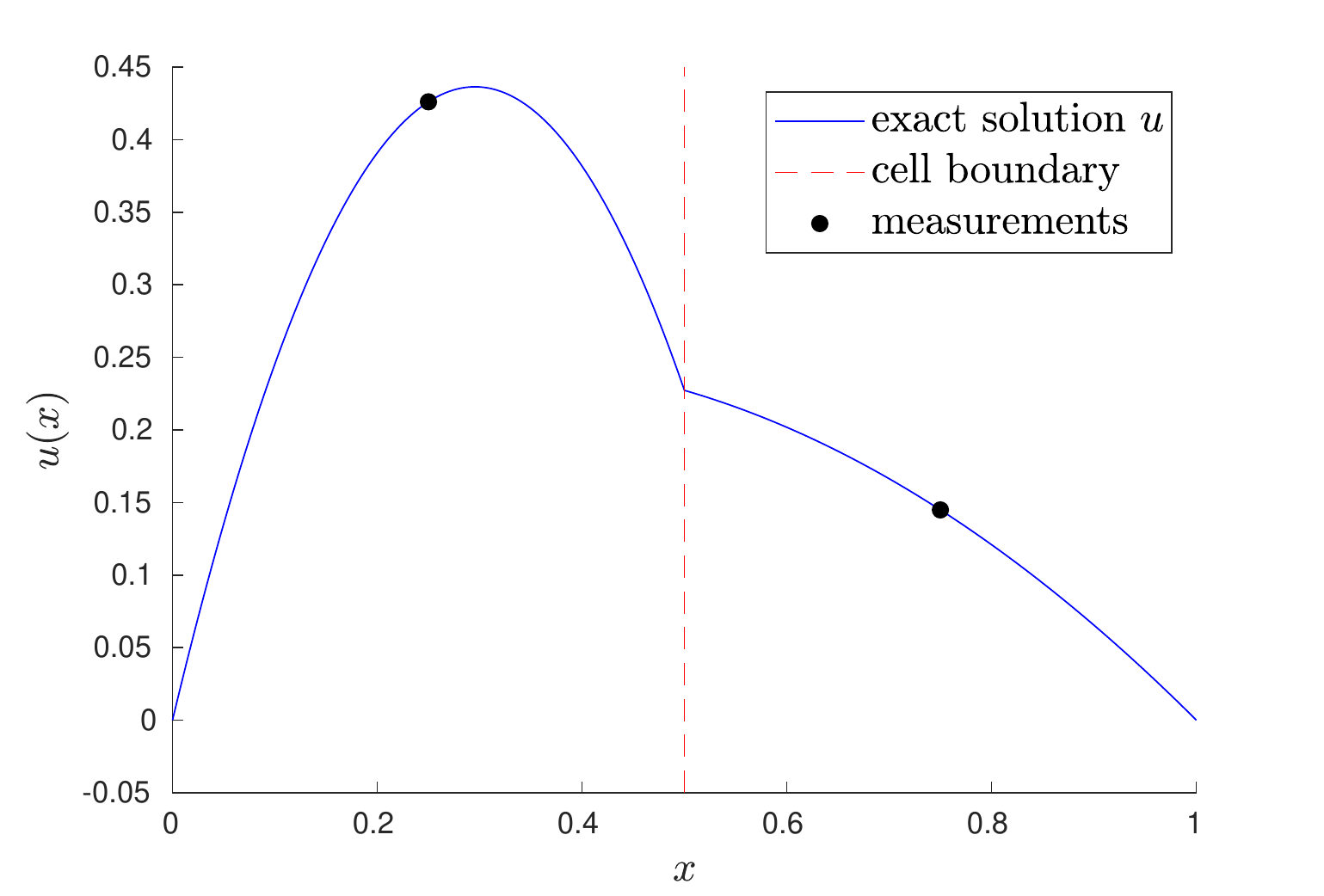}
\caption{Exact solution to 
\eqref{eq:poiss1D_1}--\eqref{eq:poiss1D_2} 
when $a(x) = 0.1$ for $0<x<1/2$ 
and $a(x) = 1$ for $1/2<x<1$, 
corresponding to ``true'' 
parameters $\hat{\theta}_0 = 0.1$ 
and $\hat{\theta}_1 = 1$,
with measurements
at $x_0 = 0.25$ and 
$x_1 = 0.75$.}
\label{fig:1D_exact_u}
\end{figure}

Let us specifically consider the case with $N=2$ parameters to understand some 
qualitative features 
of the benchmark. 
In this case,
$
c = -\frac{1}{4}(3\theta_0+\theta_1)/(\theta_0+\theta_1)
$
and 
\begin{align*}
u_0(x) &= -\frac{x^2}{2\theta_0} + \frac{3\theta_0+\theta_1}{4\theta_0(\theta_0+\theta_1) }x,
&\quad
u_1(x) &= -\frac{x^2}{2\theta_1} + \frac{3\theta_0+\theta_1}{4\theta_1(\theta_0+\theta_1) }x  + \frac{\theta_1-\theta_0}{4\theta_1(\theta_0+\theta_1)}.
\end{align*}
This solution is shown in Fig.~\ref{fig:1D_exact_u}.

Let us assume that we have measurements at $x_0=0.25$ and $x_1=0.75$. Then the ``exact'' measurements we would get are%
\footnote{Unlike in Section \ref{sec:benchmark-description}, these values are computed 
using the exact solution, 
instead of a finite element approximation.} 
\begin{align*}
\hat{z}_0 &= -\frac{x_0^2}{2\theta_0} + \frac{3\theta_0+\theta_1}{4\theta_0(\theta_0+\theta_1) }x_0, &\qquad
\hat{z}_1 &= -\frac{x_1^2}{2\theta_1} + \frac{3\theta_0+\theta_1}{4\theta_1(\theta_0+\theta_1) }x_1  + \frac{\theta_1-\theta_0}{4\theta_1(\theta_0+\theta_1)}.
\end{align*}

\begin{figure}
\includegraphics[width=0.48\textwidth]
{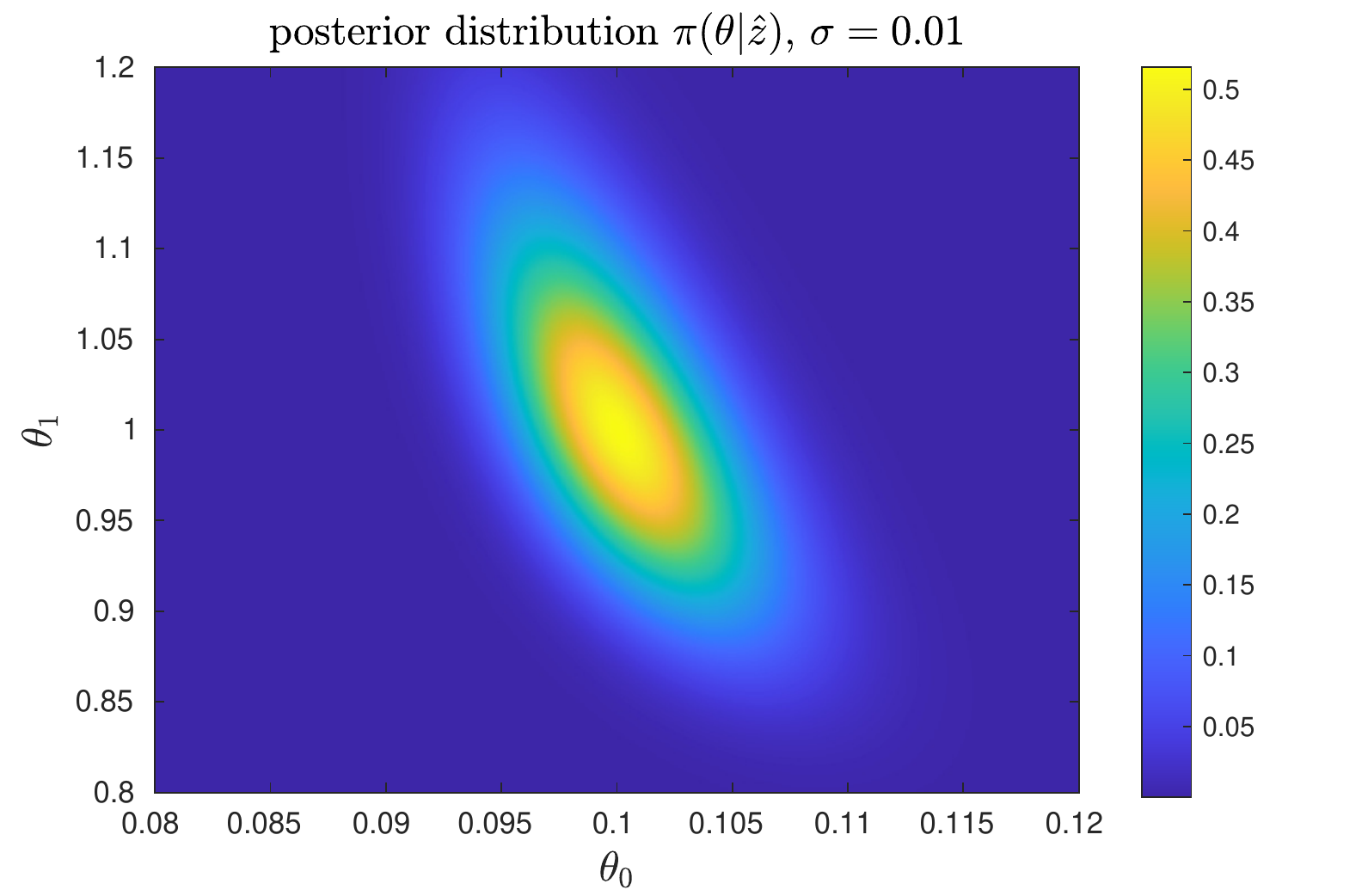}
\hfill
\includegraphics[width=0.48\textwidth]
{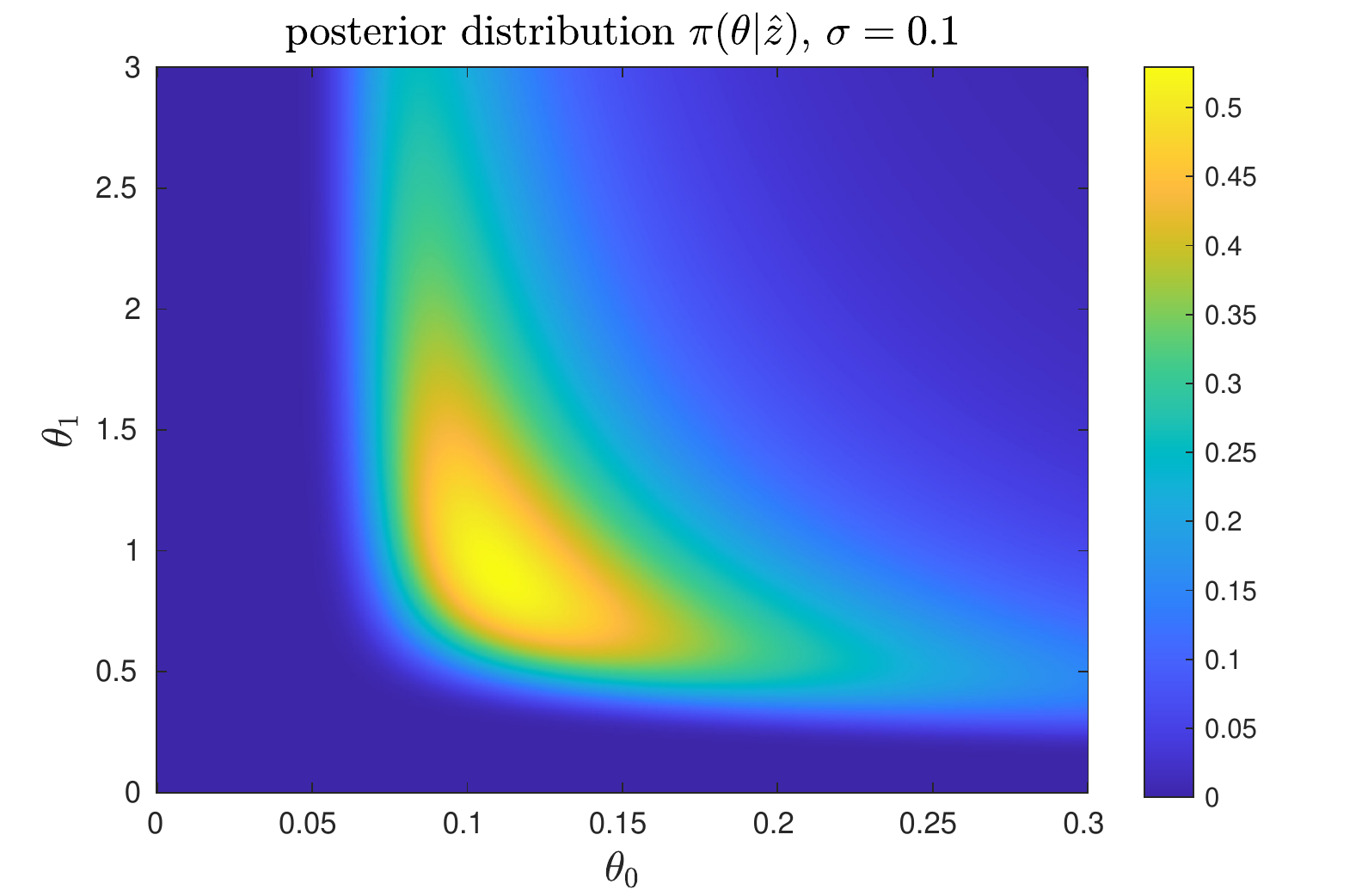}
\caption{Left: Non-normalized posterior distribution $\pi(\theta|\hat{z})$ when $\sigma = 0.01$ in the 1D model. 
Notice the distribution 
is approximately 
Gaussian for this 
very small $\sigma$.
Right: When using $\sigma = 0.1$ in the 1D model. 
For this larger value of $\sigma$, 
the posterior distribution extends to larger values of $\theta_0$ and $\theta_1$, and 
has a
non-Gaussian, ``banana'' 
shape, reminiscent of the pair histograms in Fig.~\ref{fig:pair-histogram}.}
\label{fig:1d-posterior}
\end{figure}

We can use these values to exactly and cheaply evaluate the posterior probability (without sampling), and
Fig.~\ref{fig:1d-posterior} shows $\pi(\theta|\hat{z})$ in 
this simple case, using
``true'' parameters 
$\hat{\theta}_0 = 0.1$ and 
$\hat{\theta}_1 = 1$ to define $\hat z$.
We use the prior and likelihoods defined before, 
with the 
prior standard deviation $\sigma_{\textup{pr}} = 2$ 
and two different values for the likelihood 
standard deviations used in \eqref{eq:likelihood}: 
$\sigma = 0.01$ 
and
$\sigma = 0.1$.

The figure illustrates more concisely the strong correlation between parameters that we have experimentally observed in Fig.~\ref{fig:pair-histogram}. It also illustrates that if $\sigma$ in the likelihood \eqref{eq:likelihood} is chosen small (i.e., if the measurement error is small), then the posterior is approximately Gaussian. The observation here therefore validates our choice of a relatively large $\sigma$, see also Remark~\ref{remark:sigma}.

\section{Estimating the essential sample size}
\label{sec:background}
Section~\ref{sec:information} assessed how much information is actually present in the many samples we have computed, and we have also used the results shown there in providing theoretical support for the key cost-accuracy estimate \eqref{eq:convergence}.
Many of the statistical errors described there can be estimated from the autocovariance $AC_L(s)$ defined in \eqref{eq:autocovariance}. 
The basis for this is the Markov chain central 
limit theorem. Informally, the Markov chain central limit theorem says that, for large $n$ and $N_L$, 
the running means $\langle \theta\rangle_{L,n}$ 
are approximately normally distributed with mean $\langle \theta\rangle_{\pi(\theta|\hat{z})}$ and covariance 
\begin{equation}\label{eq:en_scaling}
\textup{Cov}\left(\langle \theta\rangle_{L,n}\right) \approx \frac{\IAC}{n},
\end{equation}
where $\IAC$ is the integrated autocovariance, obtained by summing up the autocovariance $AC_L(s)$; we estimate
$\IAC$ from the data in~\eqref{eq:IAT_matrix}. Equation~\eqref{eq:en_scaling} then justifies the formula for the scaling of $e(n)^2$ below~\eqref{eq:convergence}.

To establish the effective sample size formula~\eqref{eq:ESS_perchain}, we cite the ``Delta method'' in statistics~\cite{wolter2007introduction}.
Informally, the Delta method states that, for large $n$ and $N_L$, 
a continuously differentiable function $f$ of the running means is nearly normally distributed with 
variance 
\begin{equation}\label{eq:compare_variance}
    \textup{Var}(f(\langle \theta\rangle_{L,n})) \approx  \frac{1}{n}\nabla f(\langle \theta\rangle_{\pi(\theta|\hat{z})})^T\cdot   \IAC\cdot \nabla f(\langle \theta\rangle_{\pi(\theta|\hat{z})}),
\end{equation}
provided $\nabla f(\langle \theta\rangle_{\pi(\theta|\hat{z})})$ is nonzero.
The formula~\eqref{eq:compare_variance} is obtained by Taylor expanding $f$ and applying the Markov chain central limit theorem.
For a Markov chain in which successive samples are all independent, $\IAC$ is the covariance matrix, which we estimate from the data in~\eqref{eq:covariance} and denote by $C$. Using a standard 
result on generalized Rayleigh quotients, 
\begin{align*}
    \frac{\nabla f(\langle \theta\rangle_{\pi(\theta|\hat{z})})^T\cdot   \IAC\cdot \nabla f(\langle \theta\rangle_{\pi(\theta|\hat{z})})}{\nabla f(\langle \theta\rangle_{\pi(\theta|\hat{z})})^T\cdot   C\cdot \nabla f(\langle \theta\rangle_{\pi(\theta|\hat{z})})} \le \max_{x \ne 0} \frac{x^T \cdot \IAC \cdot x}{x^T \cdot C\cdot x} = \lambda_{\max}(C^{-1}\cdot \IAC),
\end{align*}
where $\lambda_{\max}$ denotes the largest eigenvalue.
This means that the variance~\eqref{eq:compare_variance} is at most $\lambda_{\max}(C^{-1}\cdot \IAC)$ times 
the value it would take if all of the 
samples on chain $L$ were independent. 
In other words, the minimum number of 
``effectively independent samples''  is approximately $N_L/\lambda_{\max}(C^{-1}\cdot \IAC)$ for each of our chains of length $N_L$.

It is quite standard to quantify statistical 
error in Markov chain Monte Carlo simulations by using the integrated autocovariance~\cite{sokal1997monte,geyer1992practical}. In the literature, it is also common to 
determine an effective sample size by 
checking to see where autocovariances cross 
a certain threshold -- such as a small 
percentage of its initial value -- as we have 
done in Figure~\ref{fig:autocorrelation} above, where we used a threshold 
of $1\%$.

\bibliographystyle{siam}
\bibliography{paper}

\end{document}